\documentclass[a4paper,12pt]{extarticle}
\usepackage[utf8]{inputenc}
\usepackage{tikz}

\tikzset{
  mynode/.style={fill,circle,inner sep=2pt,outer sep=0pt}
}
\usepackage{enumitem}
\usepackage[english]{babel}
\usepackage{amsfonts}
\usepackage{amsthm}
\usepackage{amsmath}
\newtheorem*{theorem}{Theorem}
\usepackage{tikz}
\usepackage{graphicx}
\usepackage{tikz}
\usepackage{float}
\usetikzlibrary{matrix,positioning}
\usetikzlibrary{decorations.pathreplacing}
\usepackage{enumitem}
\usepackage{amsmath}
\usepackage[english]{babel}
\usepackage{amsthm}
\usepackage{graphicx}
\usepackage{amsthm}
\theoremstyle{remark}
\newtheorem*{lemma}{Lemma}

\title{Minimal Genus Seifert surface of $11$  crossing alternating knots}
\author{Neetal Neel}
\date{}

\begin{document}

\maketitle

\section*{Abstract}

{\scriptsize{
Kakimizu complexes have been found for several classes of links. O.Kakimizu found the Kakimizu complexes of knots with crossing number less than or equal to 10. Hatcher and Thurston found the 0-skeleton of the Kakimizu complex of 2-bridge links. Sakuma later generalized the result for special arborescent links and found the Kakimizu complexes for the same. Jessica Banks gave a complete proof of results announced by Hirasawa and Sakuma, describing explicitly the Kakimizu complexes of non-split, prime special alternating links. The Kakimizu complexes of prime, non-split alternating links have finite number of vertices. In this paper we compute the Kakimizu complexes for all $11$ crossing prime, alternating knots,  explicitly describing each of them. For most knots and links we use known algorithms. The rest of the Kakimizu complexes were found by using Murasugi sums and sutured manifold theory developed by Gabai, Scharlemann, Kakimizu and others. }}

\normalsize

\section{Introduction}

Let $K$(or $L$) be an oriented knot (or link) in $S^3.$ A Seifert surface on $K$ (or $L$) is a compact, connected, orientable surface with $K$(or $L$) as the boundary of the surface compatible with the orientation of $K$(or $L$) and which does not contain any closed component. Let $E(L)$ be the link complement defined by $S^3 - (N(L))^{\circ}$, where $N(L)$ is a regular neighborhood of $L$ homeomorphic to $L\times D^2$. $L$ will be used to denote a link (with one component links considered as knots). $K$ will be used to specifically denote a knot. 

We will call $S$ to be a Seifert surface on $L$ to denote $S$ a Seifert surface with no closed component and $L$ as the boundary of $S$ contained in $S^3$. We will abuse notation to denote $S$, a Seifert surface on $L$ lying on $E(L),$ i.e $S\cap E(L)$ with $\partial(S\cap E(L))$ is a parallel copy of $L$ lying in the boundary tori of the regular neighborhood of the link $L$, i.e   $\partial N(L).$ It will be clear from the context otherwise will be elaborated.

In 1992, Kakimizu introduced Kakimizu complex of an oriented knot. Kakimizu complex aims to explore the configuration of Seifert surfaces on links $L$ in the link complement, $E(L).$ More specifically, it aims to answer how disjoint representatives of the set of  isotopy classes of incompressible (or minimal genus) Seifert surfaces on a link $L$ is situated in the knot/link complement, expressed with the help of a complex. 

In this paper, we will consider oriented, non-split links, since $E(L)$ for non-split links are irreducible. Moreover we consider minimal genus Seifert surfaces $S$ on a link $L$ so that $S$ is an incompressible surface. Note that incompressible surfaces may not be a minimal genus Seifert surface but minimal genus Seifert surfaces are incompressible. We consider incompressible surfaces because given any surface, adding a compressible handle increases the genus and represents a new surface on the given link. This would create an infinite family of surfaces, hence an infinite complex on the given link $L$. 

Kakimizu has computed the Kakimizu complex of all knots upto $10$ crossings. In 2018, Hass, Tsvietkova and Thompson announced the result, that the number of isotopy classes of Seifert surfaces of a given genus is bounded by a polynomial depending only on the genus of the Seifert surface. This shows that the Kakimizu complex is a finite complex for alternating links. Cromwell in 1989 showed that every homogenenous link $L$ (a superset of the set of alternating links) is a $*-$product of special alternating links, $L_i$. This is equivalent to, given a reduced alternating diagram $D$ of a homogeneous link $L$ and applying Seifert's algorithm on $D,$ we obtain a Seifert surface $S$ which is the Murasugi sum of Seifert surfaces $S_i$ on $L_i.$ This served as a motivation to try to find an algorithm for computing the Kakimizu complex of homogeneous links $L.$ Now this problem turns out to be very difficult owing to the fact, that the sutured manifold theory is developed for plumbing but is absent for a general Murasugi $2n-$gon disk. Moreover the Kakimizu complex depends on the embedding site of the Murasugi disks in the corresponding Seifert surfaces on the special alternating links $L_i,$ In this paper, we compute the Kakimizu complex of all alternating $11$ crossing knots using the methods known for special classes of links. Most of the links either fall in one of the classes or could be found using Theorems of Kakimizu and Gabai.

There are a few 11 crossing links whose Kakimizu complex had to be found using sutured manifold theory. I found a new type of surface, called the parallel surface, which is a priori not a plumbing. 

I would like to thank Jennifer Schultens, my advisor,  for the constant and unwavering help throughout the work. I would also like to thank Andrew Alameda for the valuable conversations we had. 

\newpage

\section{Preliminaries}

\subsection{Sutured Manifold Theory}

This section is mostly taken from the paper by Kakimizu, \emph{Classification of the incompressible spanning surfaces for prime knots of 10 or less crossings.} For a detailed reading please refer to \cite{10.32917/hmj/1150922486}.

A \emph{sutured manifold $(M,\gamma)$} is a compact oriented $3$-manifold $M$ together with a union of finitely many pairwise disjoint annuli $\gamma\subset \partial M.$ For each component of $\gamma,$ which is an annulus, there exists a core circle which is called a \emph{suture}. The set of sutures are denoted by $s(\gamma).$ Moreover $R(\gamma)$ defined as  $\partial M - \gamma^{\circ}$ is oriented and is coherent with respect to $s(\gamma).$ Let $R_{+}(\gamma)$ (or $R_{-}(\gamma)$) be the components of $R(\gamma)$ whose normal vectors point out of $M$ (respectively point into $M.$) 
 The \emph{knot complement} of a link $L,$ is $M = S^3 - N(L),$ a compact, oriented 3-manifold and $\partial M =  \partial N(L) \cong L \times \mathbb{S}^1.$ The sutures, $s(\gamma),$ generally are parallel to the components of the link $L.$ 

Let $(M,\gamma)$ be a sutured 3-manifold. A $\gamma$-surface $S$ is a properly embedded compact oriented surface in $M$ such that $\partial S\subset (\gamma)^{\circ}$ and is isotopic to $s(\gamma),$ with no closed component. A $\gamma$-surface $S$ is said to be \emph{parallel} if there exists an isotopy $$e: (S,\partial S) \times I\to (M,\gamma)$$ with $e_0 = \text{ identity}$ and $e_1(S)\subset R(\gamma).$ An \emph{essential} $\gamma$-surface is a $\gamma$-surface which is incompressible and not parallel to a surface in $R(\gamma).$ Note that if $S$ is a $\gamma$-surface in $(M,\gamma),$ $\partial M$ is connected and $S$ is parallel to a surface in $R(\gamma),$ then $S$ is isotopic to $R_{+}(\gamma)$ or $R_{-}(\gamma).$

A $\gamma$-isotopy of $(M,\gamma)$ is a map $h:I \times M \to M$ such that $h_0 = $identity, with the property $h_t|_{R(\gamma)} = $identity, and $h_{t}|_{\gamma} = $identity for all $0\leq t\leq 1$. We say two minimal genus Seifert surfaces on a link $L$,  $S$ and $S',$ are isotopic to each other if there is an ambient isotopy  $$e:S^3\times I\to S^3$$ such that $e_0 = $identity map and $e_1(S) = S'$ with $e_t|_{L} = $identity. If $S$ and $S'$ are isotopic then there is a $\gamma$-isotopy $h:I\times M\to M$ of $(M,\gamma) = (S^3-N(L), \partial N(L))$ such that $h_1(\tilde{S}) = \tilde{S'}$ where $\tilde{S},\tilde{S'}$ denote $S$ and $S'$ restricted to the link complement ($S^3-N(L))$ respectively. We will think of a minimal genus Seifert surface $S$ on $L$ as a $\gamma$-surface on $(M,\gamma) = (S^3-N(L),\partial N(L))$ with the boundary of the surface in $\partial N(L)\cap S$, a parallel copy of $L$ in $\partial N(L).$ 

Two major examples of sutured manifolds are product sutured manifolds and the complementary sutured manifolds. 

Let $F$ be a compact, oriented $2$-manifold. Then $(M,\gamma),$ homeomorphic to $(F\times I, \partial F\times I),$ is called a product sutured manifold. If $S$ (on $L$) is a minimal genus Seifert surface in the link complement then a regular neighborhood, $$\bigg(N(S), ( S\cap \partial N(L))\times I\bigg)\subset \bigg(S^3-N(L),\partial N(L)\bigg)$$ is a product sutured manifold since $N(S) \cong S\times I$ and $S\cap \partial N(L))\cong L.$

Let $S$ be a minimal genus Seifert surface (on $L$) in the link complement. In this paper, we are  considering only links that are oriented and non-split. Let $E(L)$ denote the link complement and $(N(S),\partial S \times I)$ be the product sutured manifold. 

A sutured manifold $(M_S,\gamma)$ defined by 
$$(M_S,\gamma) = \bigg(\overline{E(L) - N(S)}, \hspace{.2cm} \overline{ \partial E(L) - ((S\cap \partial N(L))\times I)}\bigg)$$ is called the complementary sutured manifold. Given a minimal genus spanning surface $S$ (on $L$) we look for essential surfaces in the complementary sutured manifold of $S$, $(M_{S},\gamma).$ Any essential $\gamma$-surface $S'$ in $(M_S,\gamma)$ is disjoint from $S$ and is not isotopic to $S.$  

\subsection{Product Decomposition}

Let $(M,\gamma)$ be a sutured manifold. Let $D\subset M$ be a disk properly embedded in $M$ such that $\partial D$ intersects $ s(\gamma)$ in exactly at $2$ points. $D$ is said to be a \emph{product disk} in  $(M,\gamma).$ An operation on $(M,\gamma)$ called \emph{product decomposition} yields us another sutured manifold $(M',\gamma').$ Each component of $R(\gamma')$ is incompressible if and only if every component of $R(\gamma)$ is incompressible. We denote this operation as 

$$(M,\gamma) \xrightarrow{D} (M',\gamma')$$

Let $(M,\gamma)$ be a sutured manifold. Note that the complementary sutured manifold of a non-split link is irreducible. Moreover a $3$-submanifold of an irreducible $3$-manifold is irreducible.  

Let $S\subset M$ be a $\gamma$-surface. By a $\gamma$-isotopy we can isotope $S$ such that $\partial S = s(\gamma)$ and $S$ intersects transversally $\partial D$ in an arc connecting the $2$ points of intersection of $s(\gamma)$ and $D.$ After product decomposition, cut through the arc, and we get a new $\gamma'$-surface $S_D$ in $(M',\gamma').$

This is from [\cite{10.32917/hmj/1150922486} , Lemma 1.3]

\begin{theorem}

Let $(M,\gamma)\xrightarrow{D}(M',\gamma')$ be a product decomposition. Suppose
that $M$ is irreducible and $\partial M'$ is connected. Then for each essential $\gamma$-surface
$S \subset M$, the $\gamma'$
-surface $S_D$ is essential. Moreover if $S,S'$ are $\gamma-$isotopic to each other then $S_D, S_D'$ are $\gamma'$-isotopic as well.

\end{theorem}

Let $S_{(M,\gamma)}$ be the isotopy classes of essential surfaces $S\subset  (M,\gamma).$ Two essential surfaces $S, S'$ are in the same isotopy class if they are $\gamma$-isotopic to each other.  Let $(M,\gamma)\xrightarrow{D}(M',\gamma')$ be product decomposition. Consider the map:

$$S_{(M,\gamma)}\xrightarrow{\widetilde{D}} S_{(M',\gamma')} \text{ given by } \widetilde{D} (S') = S_D'.$$ The previous theorem indicates that the above map is well defined.

This is [\cite{10.32917/hmj/1150922486} ,Lemma 1.4]

\begin{theorem}

Suppose that $M$ is irreducible and $\partial M'$ is connected. Then the map: $S_{(M,\gamma)}\xrightarrow{\widetilde{D}} S_{(M',\gamma')}$ is a bijection.
    
\end{theorem}

For fibred links $L$ with fibre $S,$  the complementary sutured manifold $(M_S,\gamma)$ is a product sutured manifold homeomorphic to $(S\times I, L\times I).$

Hence the next theorem is useful in this context:

\begin{theorem}

Let $(M,\gamma)\xrightarrow{D} (M',\gamma')$ be a product decomposition. Let $M$ be irreducible and that $(M',\gamma')$ has two components, $(M_1,\gamma_1)$ and $(M_2,\gamma_2).$ Let $(M_2,\gamma_2)$ be a product sutured manifold and $\partial M_1$ be connected. Then we have a bijection:

$$S_{(M,\gamma)}\xrightarrow{\widetilde{\bar{D}}} S_{(M_2,\gamma_2)}$$ is a bijection.
    
\end{theorem}

Finally a fundamental theorem proved by Kakimizu using Haken and Hempel and Waldhausen's theorems. \cite{667e8b29-d581-32a3-ab92-1d9363b3b401}, 

\begin{theorem}

Let $X$ be a connected Haken $3$-manifold such that $\partial X$ is a union of incompressible tori. Let $Y$ be a compact irreducible $3$-submanifold of $X$ (possibly disconnected) such that each component of $\partial Y$ is a properly embedded incompressible surface in $X$. Let $F$ and $F'$ be two properly embedded
orientable incompressible surfaces in $X$ (possibly disconnected) which satisfy the
following properties $(1)–(4)$. Then there is an isotopy $h_t$ of X keeping $Y$ fixed
so that $h_0= id$ and $h_1(F) = F'.$
\begin{enumerate}[label=(\arabic*)]
    \item $F\cup F' \subset X-Y.$

    \item Each component of $\partial X $ contains at most one component of $\partial F$ and $F$ has no closed components.

    \item There is a homotopy $f: F\times I\to X$ such that $f_0 = id$ and $f_1:F\to F'$ is a homeomorphism and $f(\partial F\times I)\subset \partial X,$

    \item There is no component of $F$ which is parallel to a component of $\partial Y.$

\end{enumerate}
    
\end{theorem}

The proof is in \cite{10.32917/hmj/1150922486} .

An application of this theorem is as follows. Let $X$ (from the theorem) be the link complement $E(L)$ for a non-split link $L.$ Then $X$ is an irreducible Haken manifold with $\partial X$ are incompressible tori with the core circles are link components. Let $S$ be an essential surface on $L.$ Consider the complementary sutured manifold $Y = (M_S,\gamma).$ Let $S', S''$ be two essential embedded $\gamma$-surfaces in $(M_S,\gamma)$ such that $S'$ and $S''$ are isotopic in $X.$ Then it satisfies all the conditions $(1) - (4),$ the theorem implies that there is a $\gamma$-sutured isotopy of $S'$ and $S''$ in $(M,\gamma).$ 

\subsection{Murasugi sums and Plumbings}

An oriented surface $S$ is said to be a Murasugi sum of $S_1$ and $S_2$ if there exists a splitting sphere $S^2\subset S^3$ splitting $S^3$ into two balls $V_1$ and $V_2$ ($V_1\cup V_2 = S^3$, $V_1\cap V_2 = S^2$ ) with $$S_1\subset V_1,S_2\subset V_2; \hspace{1cm} S_1\cup S_2 = S, S_1\cap S_2 = S\cap (S^2 = \partial V_i) = D$$ where $D$ is an embedded $2n$-gon on $S$.  The disk $D$ is said to be the $2n$-Murasugi disk for the Murasugi sum $S = S_1\cup_D S_2.$ 

Let $\partial S = L$ and $\partial S_i = L_i.$ Then $S\cap E(L) \subset E(L)$ is also said to be the Murasugi sum of $S_1\cap E(L_1)$ and $S_2\cap E(L_2).$ Also $L$ is said to be the Murasugi sum of $L_1$ and $L_2.$ Let $S^{c} = S_1\cup_{D^{c}}S_2$ with $D^c = (S^2 = \partial V_1)- D.$ This is obtained by $S^{c} = (S-D)\cup D^{c}.$ Note that $S^{c}$ is also the Murasugi sum of $S_1$ and $S_2.$ This surface is called the dual of $S.$ By an isotopy of $S^3,$ keeping the link $L$ fixed, we can isotope $S^c$ such that $S\cap S^c = \phi.$ For $n=2,$ the Murasugi disk is a rectangle and $4$-Murasugi sum is called plumbing. 

Gabai showed the following 2 properties for Murasugi sum operation. Assume $S$ on link $L$ be the Murasugi sum of $S_1$ and $S_2$ on $L_1$ and $L_2$ respectively.  

\begin{itemize}
    \item $S$ is a minimal genus Seifert surface on $L$ if and only if $S_1$ and $S_2$ are both minimal genus Seifert surfaces on $L_1$ and $L_2$ respectively.

    \item $L$ is a fibred link with $S$ being the fibre if and only if $L_1$ and $L_2$ are fibred links with $S_1$ and $S_2$ being the fibres for the respective fibred links.
\end{itemize}

\newpage

\section{Definition of Kakimizu Complex}

Let $L$ be a non-split oriented link. Seifert's algorithm applied on a reduced alternating diagram $D$ of $L,$ yields a Seifert surface $S$ on $L.$ Let $E(L)$ denote the link complement. $E(L)$ is an irreducible $3$-manifold since $L$ is non split. Consider the surface $S\cap E(L)\subset E(L).$ This is a properly embedded compact surface in $E(L)$ with the boundary, $\partial (S\cap E(L))\subset \partial E(L)$ homeomorphic to $L.$ Let $S$ be a Seifert surface in $E(L)$ which is compressible in $E(L).$ Then we can use the compression disks to compress the surface until the resultant surface is incompressible. Since every minimal genus Seifert surface on a link $L$ is incompressible, we consider minimal genus Seifert surface for the rest of the paper. This establishes that every oriented non-split link $L$ admits a minimal genus Seifert surface.

Let $S,S'$ be Seifert surfaces on a link $L.$ $S$ and $S'$ are said to be isotopic if there exists an isotopy fixing the link $L,$ taking $S$ to $S'.$ Let $e$ be the isotopy  $e:I\times S^3\to S^3$ with $e_0 = id$ ; $e_1(S) = S'$ and $e_t|_{L} = id$ for all $0\leq t\leq 1.$ Given an isotopy $e,$ of $S^3$  taking $S$ to $S'$ with the link $L$ fixed, there exists an isotopy $e'$ of $E(L)$ such that the surface $S\cap E(L)$ is taken to  $S'\cap E(L)$ with $e'_t(\partial E(L)) = \partial E(L)$ for all $0\leq t\leq 1.$  

The $0$-skeleton of the Kakimizu complex of a link $L$ are isotopy classes of minimal genus Seifert surfaces on the link $L.$ The $1$-skeleton of the Kakimizu complex of $L$ are edges between two vertices $v$ and $v'$ if there exists representative surfaces $S$ and $S'$ respectively such that $S$ and $S'$ are disjoint in $E(L).$ Kakimizu complex is a flag complex. If there are $n+1$ vertices such that every pair of vertex has an edge then there is an $n$ simplex on the corresponding $n+1$ vertices.

\newpage

\section{Computing the Kakimizu Complex of knots}

\subsection{Kakimizu complex of fibred knots}

In 1972, W.Whitten showed in \cite{WHITTEN1973373} that for a fibred link $L$ there is a unique incompressible surface for a fibred link $L$, unique in the sense of the definition of Kakimizu complex. Hence $MS(L) = \{[S]\}$ which implies that the Kakimizu complex of a fibred knot is a point.

We derive the set of fibred knots from KnotInfo \cite{knotinfo}.

There is an algorithm to check if a homogeneous link $L$ (that includes the class of alternating links) is fibred or not. In 1989, Cromwell \cite{https://doi.org/10.1112/jlms/s2-39.3.535} showed that a homogeneous link $L$ is a *-product of special alternating links, $L_i$. Given an oriented diagram $D$ of a homogeneous link $L$, applying Seifert's algorithm to the diagram $D$ yields a Seifert surface $S$ which is a Murasugi sum of Seifert surfaces $S_i$ on special alternating links $L_i$. Gabai \cite{Gabai1983TheMS} showed that the Murasugi sum of two links, $L_1$ and $L_2$ is fibred if and only if both $L_1,$ and $L_2$ are fibred.

Each homogeneous link $L$ is a Murasugi sum of special alternating link,  Jessica Banks \cite{banks2012minimal} has given an algorithm to detect fibredness of a special alternating link: Let $D$ be a reduced diagram of a special alternating link $L.$

Let $S(D)$ be the partition of $S^2$ into black and white regions by applying Seifert's algorithm and each Seifert disk constitutes a black region in $S^2$. Let  $G(D)$ be the planar graph with a vertex in each white region and an edge through each crossing. A special alternating link $L$ is fibred if given a reduced diagram $D$ of $L,$ we can reduce the graph $G(D)$ to a single vertex using the following moves:

\begin{itemize}
    \item Delete a loop.

    \item Contract an edge if one of the endpoint is a vertex of valence 2. 
\end{itemize}

If every summand of a Murasugi sum of links is fibred then the Murasugi sum of links is a fibred link.

\clearpage

{\bf{The list of fibred links for $11$ crossing alternating links is as follows:}}

 $$11_3,11_5,11_7,11_9,11_{14},11_{15},11_{17},11_{19},11_{22},11_{24},11_{25},11_{26},11_{28},11_{33},11_{34},11_{35},11_{40},11_{42},$$ 
 $$11_{44},11_{47},11_{51},11_{53},11_{55},11_{57},11_{62},11_{66},11_{68},11_{71},11_{72},11_{73},11_{74},11_{76},11_{79},11_{80},11_{81},$$
 $$11_{82},11_{83},11_{86},11_{88},11_{92},11_{96},11_{99},11_{106},11_{108},11_{109},11_{112},11_{113},11_{121},11_{125},11_{126},11_{127},$$
 $$11_{128},11_{129},11_{131},11_{139},11_{142},11_{146},11_{147},11_{151},11_{156},11_{157},11_{158},11_{159},11_{160},11_{162},11_{163},$$
 $$11_{164},11_{170},11_{171},11_{174},11_{175},11_{176},11_{177},11_{179},11_{180},11_{182},11_{184},11_{189},11_{194},11_{196},11_{203},$$
 $$11_{206},11_{209},11_{215},11_{216},11_{217},11_{218},11_{221},11_{223},11_{228},11_{231},11_{232},11_{233},11_{239},11_{248},11_{250},$$
 $$11_{251},11_{252},11_{253},11_{254},11_{255},11_{257},11_{259},11_{261},11_{264},11_{266},11_{267},11_{268},11_{269},11_{274},11_{277},$$
 $$11_{281},11_{282},11_{284},11_{286},11_{287},11_{288},11_{289},11_{293},11_{300},11_{301},11_{302},11_{305},11_{306},11_{308},11_{314},$$
 $$11_{315},11_{316},11_{326},11_{330},11_{332},11_{346},11_{348},11_{350},11_{351},11_{367}.$$

\vspace{1cm}

{\bf{The Kakimizu complex of fibred links is a single vertex.}}

    \begin{tikzpicture}
    \draw[olive,thick,latex-latex] (0,0) 
    node[pos=0,mynode,fill=red,label=above:\textcolor{red}{$T$}]{};
  \end{tikzpicture}

\subsection{Kakimizu Complex of special alternating knots}

In 2012, Jessica Banks \cite{banks2012minimal} gave an explicit algorithm to find the Kakimizu complex of a prime, non-split, oriented, special alternating link $L$. A special alternating link $L$ admits a reduced, oriented, alternating diagram $D$, with the property that every Seifert circle is innermost in itself, when Seifert's algorithm is applied on $D.$

In 1991, Thistlewait and Menasco \cite{6129636b-b747-3d49-9531-0b1941964c84} proved the flyping conjecture.  \emph{Given any two reduced alternating diagrams $D_1$ and $D_2$ of an
oriented, prime alternating link $L$ then $D_1$ can be transformed to $D_2$ by applying a sequence of flypes.}

\emph{Every Seifert surface on a non-split, prime, oriented, special alternating link $L$ is isotopic to a Seifert surface obtained by applying Seifert's algorithm on a reduced, oriented, prime special alternating diagram $D$ of L.}

Hirasawa-Sakuma \cite{doi:10.1142/9789814529891}, in 1990 announced the algorithm to find the Kakimizu complex for a special alternating link. Jessica Banks \cite{banks2012minimal} and Joshua Greene \cite{article} proved the result about the characterization of Seifert surfaces on a special alternating knot independently.  

To compute the Kakimizu complex of a link $L,$ it suffices to compute the $0$ and the $1$ skeleton of the complex. Since the Kakimizu complex is a flag complex, the $1$ skeleton completely describes the Kakimizu complex of the link $L$. 

\vspace{.5cm}

\emph{Kakimizu complex of a non-split, oriented link is connected.} Proved independently by Thompson and Scharlemann \cite{10.1112/blms/20.1.61} in 1988 and Kakimizu \cite{10.32917/hmj/1150922486} in 1990, It led to the  conjecture that the Kakimizu complex of a non-split, prime link $L$ is contractible.  This has been shown to be true for different classes of links with Jessica Banks \cite{banks2012minimal} proving the result announced by Hirasawa-Sakuma for special alternating links.This was finally proved by Przytycki and Schultens \cite{25203298-73d8-3f83-8061-a5e5d88da654} in 2010. 

The Kakimizu complex of special alternating links is given explicitly by \emph{Jessica Banks} in \cite{banks2012minimal}. We explain the algorithm here briefly.

Since the Kakimizu complex of a link $L$ is connected, it suffices to find the largest simplices containing a generic vertex $v$ of the Kakimizu complex of a special alternating link $L$.

Let $L$ be a special alternating link. Let $v$ be a generic vertex of the Kakimizu complex of $L$ and let $S$ be a representative surface (in it's isotopy class) of $v.$ Since $L$ is a special alternating link, $S$ is isotopic to a surface which can be realized by applying Seifert's algorithm on a reduced, non-split, prime, oriented, alternating diagram $D.$

Any generic vertex $v$ of a non-split, special, alternating link $L$ is given by prescribing a reduced, alternating, oriented diagram $D$ of the link $L$ and applying Seifert's algorithm on $D$ to obtain a representative surface $S.$

The algorithm to find the largest complex containing the vertex $v$ is described as follows. Let $L$ be a prime, special alternating link $L.$  Given $v,$ let $D$ be the reduced alternating diagram of $L$ such that applying Seifert's algorithm on $D$ yields a representative Seifert surface $S$ of $v.$ 

\begin{itemize}
    \item Consider the Seifert graph associated to the oriented diagram $D.$ $\theta$-graph is obtained from the Seifert graph by the following process:
    
    \begin{itemize}
        \item Give a weight 1 to every edge in the Seifert graph.

        \item Identify the pair of edges with same end vertices such that they bound a bigon region in between them in $D.$ [We consider the diagram $D$ to be embedded in $S^2$]. Replace each bigon region with one edge between the two vertex, with the new edge weight being the sum of the two weights of the previous two edges. Continue till we get rid of all bigon regions.

        \item Add edges to a pair of vertices in the resulting graph with weight 0 only if there is another edge between the two vertices and the addition of this new edge doesn't create a bigon region in the graph.
        
        \item Consider the sub graph of the resulting graph that contains edges such that the endpoint vertices bound multiple edges. 
        
        \item We call this weighted graph, the $\theta$-graph of $S.$

    \end{itemize}

     \item A $\theta$- graph of a diagram $D$ of $L,$ embedded in $S^2$ divides $S^2$ into regions. The boundary of a region $r$ has edges of the $\theta$-graph. We can assign a signature to each boundary edge depending on the orientation of the diagram $D.$

   \item We fix an order of the edges of the $\theta$-graph and denote the surface $S$ as a tuple $(e_1,e_2,\dots e_n)$ with $e_i$ being the weight on the $i$th edge according to the weight.

   \item To get the maximal simplices containing the generic vertex $v,$ we need to apply all the regions of the $\theta$-graph. 
   
   \item Applying a region to the $\theta$-graph yields a non-isotopic Seifert surface which is disjoint from $S.$ Applying a region amounts to applying all the flypes associated to the weights and signature of the boundary edges of the region. 

   \item Let $r$ be a region. Let $(w_1,w_2,\dots,w_k)$ be the boundary edges of the region with the signature of the edges(w.r.t the region $r$) be $(s_1,s_2,\dots,s_k); $ $s_i = \pm 1.$

   \item After applying the region $r,$ we get a new surface, with edge weights being $w_i\to w_i +s_i$ and for non-boundary edges $e_i\neq w_j, $ $e_i\to e_i.$

   \item Apply all the regions (each region once). This implies we have applied every flype twice hence we will get back the same surface we started with.

   \item A cycle (with every region of the $\theta$-graph applied) constitutes a largest simplex containing $v.$

\end{itemize}

Since the Kakimizu complex is connected, we can obtain the Kakimizu complex of a link $L$ by choosing a generic vertex $v$, computing the largest simplices for each vertex (there is finitely many vertices in the Kakimizu complex) and we obtain the Kakimizu complex of the link $L$.  

$K = 11_{237}$ is a special alternating knot and the Kakimizu complex is computed here to demonstrate the algorithm. $(0,1,0)$ denotes a generic vertex in the Kakimizu complex. In this example, by symmetry if we compute the largest simplices containing a vertex, we will obtain the same $2$-simplex. The Kakimizu complex of $K$ is the $2$-simplex, containing the $3$ vertices $(0,1,0);(0,0,1);(1,0,0).$

\begin{figure}[H]
  \includegraphics[scale=.25]{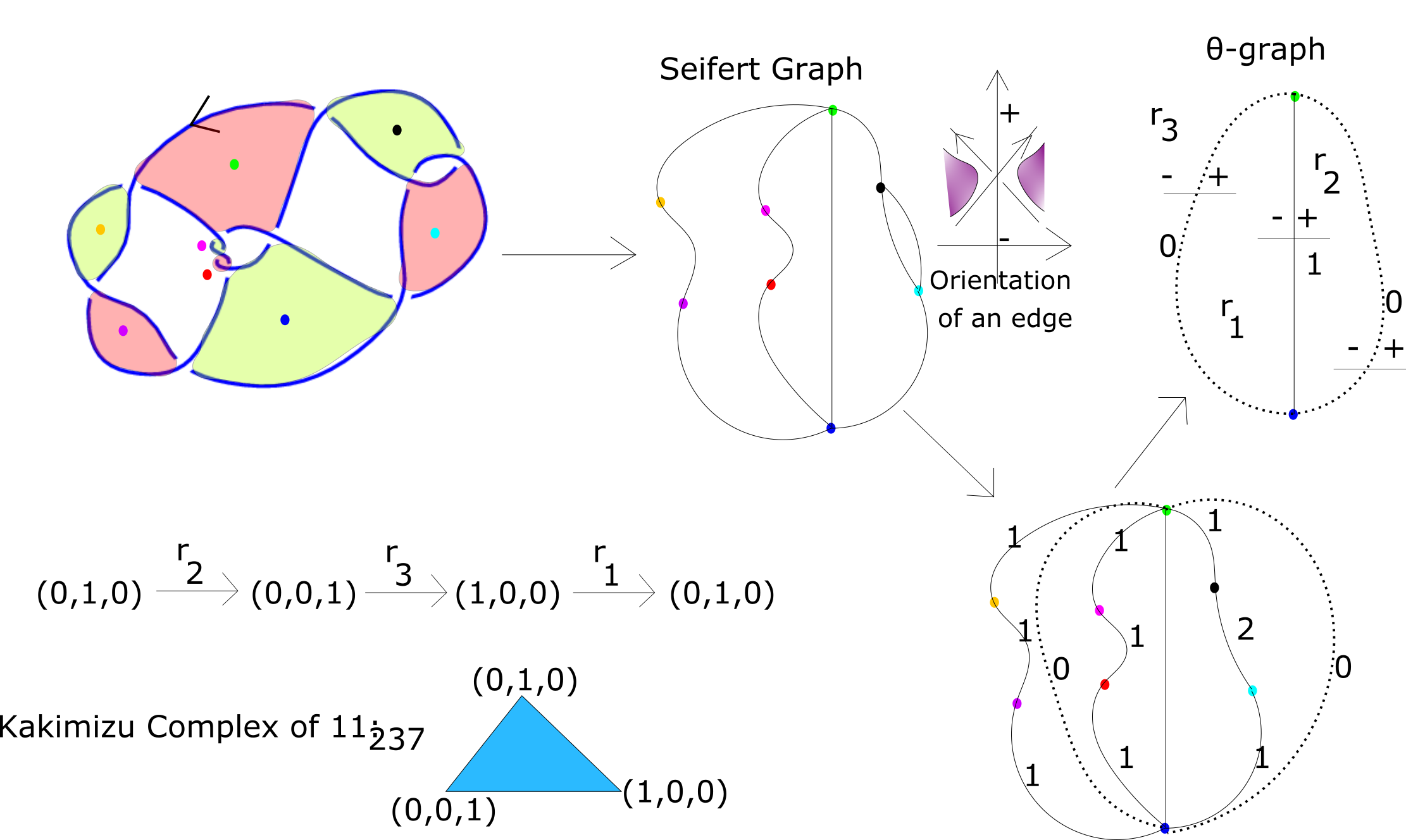}
  \caption{The Kakimizu complex of $11_{237}.$ The knot $11_{237}$ is a special alternating knot.}
  \label{}
\end{figure}
    
The list of special alternating links for $11$ crossing alternating links (with the Kakimizu complex) is as follows:

\begin{itemize}
    \item 
    $$K = 11_{43},11_{123},11_{124},11_{200},11_{227},11_{240},11_{241},11_{244},11_{245},11_{263},11_{291},11_{292}, $$
    $$11_{298},11_{299},11_{318},11_{319},11_{320},11_{329},11_{338},11_{354},11_{361}$$

are spanned by a unique minimal genus Seifert surface up to isotopy.

The Kakimizu complex of the knots are a single vertex.

    \begin{tikzpicture}
    \draw[olive,thick,latex-latex] (0,0) 
    node[pos=0,mynode,fill=red,label=above:\textcolor{red}{$T$}]{};
  \end{tikzpicture}

\item $K = 11_{94}.$ The $\theta$-graph contains 2 regions. The Kakimizu complex is  

\begin{tikzpicture}
    \draw[olive,thick,latex-latex] (0,0) -- (1,0)
    node[pos=0,mynode,fill=red,label=above:\textcolor{red}{$T_1$}]{}
    node[pos=1,mynode,fill=blue,text=blue,label=above:\textcolor{blue}{$T_2$}]{};
  \end{tikzpicture}

\item $K = 11_{237}.$ The $\theta$-graph contains $3$ regions.

\begin{tikzpicture}

\fill[fill=black!20] (0,0)--(1,0)--(1,1);

    \draw[olive,thick,latex-latex] (1,1) -- (1,0)
    node[pos=0,mynode,fill=green,label=above:\textcolor{green}{$T_3$}]{}
    node[pos=1,mynode,fill=blue,text=blue,label=above:\textcolor{blue}{$T_2$}]{};

    \draw[olive,thick,latex-latex] (0,0) -- (1,0)
    node[pos=0,mynode,fill=red,label=above:\textcolor{red}{$T_1$}]{}
    node[pos=1,mynode,fill=blue,text=blue,label=above:\textcolor{blue}{$T_2$}]{};

  \draw[olive,thick,latex-latex] (0,0) -- (1,1)
    node[pos=0,mynode,fill=red,label=above:\textcolor{red}{$T_1$}]{}
    node[pos=1,mynode,fill=green,text=blue,label=above:\textcolor{green}{$T_3$}]{};
  \end{tikzpicture}

\item $K = 11_{340}.$ The $\theta$-graph contains $2$ regions. The Kakimizu complex is:

\begin{tikzpicture}
    \draw[olive,thick,latex-latex] (0,0) -- (1,0)
    node[pos=0,mynode,fill=red,label=above:\textcolor{red}{$T_1$}]{}
    node[pos=1,mynode,fill=blue,text=blue,label=above:\textcolor{blue}{$T_2$}]{};
  \end{tikzpicture}

\end{itemize}

\subsection{Kakimizu Complex of 2 bridge knots}

Let $L$ be a $2$ bridge knot. Reidemeister \cite{Reidemeister1935HomotopieringeUL}, Schubert and Seifert \cite{Schubert1956} classified all $2$ bridge knots. Every $2$ bridge knot is associated to a reduced rational number $p/q.$ Two $2$ bridge knots with the associated bridge index $p/q$ and $p'/q'$ are equivalent if and only if $q=q'$ and $p= p^{\pm 1} (mod$ $q).$ For our purpose, given an index for a 2-bridge knot, 

\[
   \text{ Given } 0\leq p/q<1 \text{ we consider }: 
\begin{cases}
    p/q-1,& \text{if } p,q \text{ both odd }\\
    p/q,              & \text{otherwise}
\end{cases}
\]

They represent the same respective knot but the rational number $p/q$ has at least one of the entry (p,q) is even.

\emph{Every rational number $p/q$ with one even entry p or q, yields a continued fraction expansion with all even integers.}

Example:

Let $K = 11_{13}$ is a 2-bridge knot with the index $28/61.$

The continued fraction expansion of $28/61$ is given by:

\[
\mathrm{28/61}=
\cfrac{1}{2-\cfrac{1}{-6-\cfrac{1}{-2-\cfrac{1}{2}}}}
\]

Hence according to the continued fraction expansion we can denote 

$$28/61 = [2,-6,-2,2].$$

The knot could be isotoped to a 2-bridge diagram (with respect to the height function) such that it is ordered as

The $MS(L)$ (the 0-skeleton) of the Kakimizu complex can be computed from the results of Hatcher and Thurston \cite{Hatcher1985} . They showed that every minimal genus Seifert surface on a $2$ bridge link is given by the inner or outer plumbing of full twisted $n$ bands, the numbers determined from the continued fraction expansion of the 2-bridge index of the knot.

Example: 

Any surface on the knot $K = 11_{13}$ with bridge index $28/61$ is given by 

$S = B_2 \cup_{D_1} B_6\cup_{D_2} B_2\cup_{D_3} B_2$ with $B_i$ being bands with $i$ twists and $D_i$ are the plumbing disks. Every minimal genus Seifert surface on $K$ is given by inner or outer plumbings on the 3 plumbing disks. In our example since $B_2$ is the Hopf band, it's fibred, and $B_6$ is the unique incompressible, non-fibred minimal genus Seifert surface on the $3$ full twisted link, the link is not fibred and since $IS(\partial B_6) = [B_6],$ hence deplumbing each Hopf band mounts to product decomposition, and invoking the bijection from Boileau and Gabai's result with $S = S_1\cup_D F$, where $\partial F$ is a fibred link with fibre $F$, [$L_i =\partial S_i$]

$$\phi:IS(L,S) \to IS(L_1,S_1 = B_6\cup_{D_2} B_2\cup_{D_3} B_2) \to IS(L_2,S_2 = B_6\cup_{D_2} B_2) $$

$$\to IS(L_3,S_3 = B_6)  =\Phi.$$ This implies $IS(L,S) = \phi$ so $S$ is the unique spanning surface on $L.$

Sakuma[] has computed the Kakimizu complex of special arborescent links, generalizing the Kakimizu complex of $2$ bridge links. We are going to give the algorithm to find the Kakimizu complex for $2$ bridge links.

The $0-$ skeleton of a $2$ bridge link has already been found by Hatcher and Thurston. To find the Kakimizu complex of the link $L,$ we will find the largest complex containing a generic surface $S.$

Let $L$ be a $2$ bridge link with a bridge index $p/q.$ Let us consider the continued fraction of $p/q (p/q-1, \text{ if p,q both odd})$ such that each entry is even. Let $p/q = [e_1,e_2,\dots,e_k],$ with each $e_i$ being even.

To find $MS(L),$ the 0-skeleton of the Kakimizu complex of a $2-$ bridge link, denote $0$ as the inner plumbing and $1$ as the outer plumbing.

Any minimal genus Seifert surface $S = S_1\cup_{D_1}S_2\cup_{D_2}\dots\cup_{D_{k-1}}S_k$ with $S_i= B_{e_i}$ are the $e_i$ twisted bands and $D_i$ being the plumbing disk, either inner or outer with respect to a standard 2-bridge diagram.

$B_{2n}$ is a band on the $2n$ twisted link. $n=1$ represents the Hopf link with $B_2$ being the Hopf band. The Hopf link is fibred with fibre $B_2.$ For $n\geq 2,$ $B_{2n}$ is a unique minimal genus Seifert surface on the $2n$ twisted link.

Each surface in $MS(L)$ can be represented by a $k-1$ tuple of $0$ or $1.$

For example, let $L = [4,6,4].$ Then we have 2 plumbing disk, hence each  surface is a 2 tuple. The surfaces are denoted by,

$$(0,0) ; (0,1) ; (1,0); (1,1).$$

In case there is a Hopf band at the end of the chain, then  $$(0,\dots) = (1,\dots) (\text{ if Hopf band is attached at the start of the chain})$$ or $$(\dots,0) = (\dots,1) (\text{ if Hopf band is attached at the end of the chain}).$$

If the Hopf band is attached in the middle between the kth and $k+1$th plumbing disk, 

$$(\dots,s_k,s_{k+1},\dots)$$ then

$$(\dots,0,0,\dots) = (\dots,1,1,\dots).$$

This relationship holds for every Hopf band in the chain.

The algorithm to find the maximal simplex containing a generic vertex $v$ is given:

\begin{itemize}
    \item For a $2$-bridge link, let $p/q$ $(p/q-1$ for p,q both odd) be the 2-bridge index. Let $L = [e_1,e_2,\dots, e_n]$ be the even continued fraction expansion of $L.$ Let $v$ be a generic vertex of the Kakimizu complex and let $S$ be the surface denoted by $v 
 = (s_1,s_2,\dots,s_{n-1}).$ 

\item To find a maximal simplex, we need to start from $v$ and apply the \emph{component surfaces} to find a maximal cycle. Each such maximal cycle represents the maximal simplices containing $v.$

\item Let $L = [e_1,e_2,\dots,e_n],$ with $v = (s_1,s_2,\dots s_{n-1}).$ This implies $$S = B_{e_1}\cup_{D_1} B_{e_2}\cup_{D_2}B_{e_3}\dots\cup_{D_{n-1}} B_{e_n}.$$ along with the fact that $D_i$ is an inner plumbing if $s_i = 0$ and $D_i$ is an outer plumbing if $s_i = 1.$ 

\item Given a surface $S,$ denoted by $(s_1,s_2,\dots,s_{n-1})$, we can apply a component surface $B_{e_k}$ $(1<k<n)$ only if $s_{k} = s_{k+1}.$ If the condition is satisfied then applying $B_{e_k}$ on $S$ yields us $S'$ denoted by $$(s_1,\dots,s_{k-1}, (s_{k}+1) (\text{ mod } 2 ), (s_{k+1}+1) (\text{ mod } 2 ), s_{k+2},\dots s_n).$$

\item The two end component surfaces could be applied at any instant. If we apply $B_{e_1}$, on a surface $S,$ denoted by $(s_1,s_2,\dots,s_{n-1}),$ it yields us $S'$ denoted by

$$(s_1+1(\text{ mod } 2),s_2,\dots,s_{n-1}).$$

Similarly, if we apply $B_{e_n}$, on a surface $S,$ denoted by $(s_1,s_2,\dots,s_{n-1}),$ it yields us $S'$ denoted by

$$(s_1,s_2,\dots,s_{n-2},(s_{n-1}+1) (\text{ mod } 2)).$$

\item Given a generic vertex $v,$ let $S$ be the surface denoted by $(s_1,s_2,\dots,s_n).$

Consider a chain of surfaces, starting from $S$ where each component surface $B_{e_k},$ $1\leq k\leq n$ is applied exactly once (in some order). Note that  the end of the chain should exactly be $S.$ This is because any entry $s_i,$ is affected when $B_{e_i}$ and $B_{e_{i+1}}$ were applied. For any other $B_{e_{k}},$ $s_i$ remains unaffected. Applying twice, we get $$\bigg((s_i + 1) (\text{ mod }2) + 1\bigg) (\text { mod } 2) = s_i.$$ 

Hence we get a maximal cycle of surfaces starting from $S,$ containing $n$ distinct surfaces. The $n$-vertices in the Kakimizu complex span a maximal simplex containing $v.$

\item The Kakimizu complex of a non-split link is connected. Given a generic vertex $v,$ we are able to find the maximal simplices containing $v.$ Since there are finitely many vertex in the $0$-skeleton in the Kakimizu complex, applying this on every vertex yields us the Kakimizu complex.  
    
\end{itemize}

An example to illustrate the algorithm for $2$-bridge knots

\begin{figure}[H]
  \includegraphics[scale=.14]{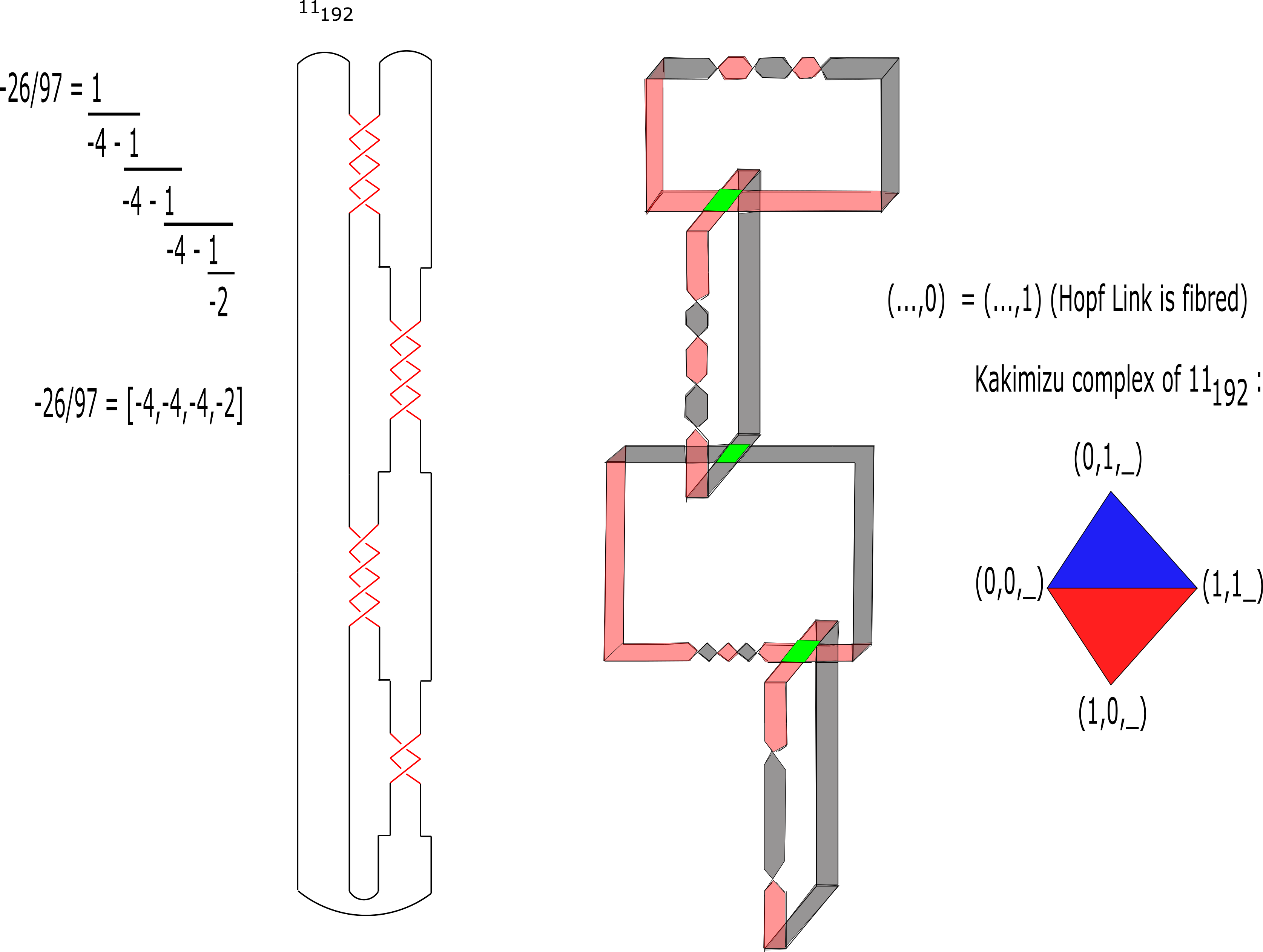}
  \caption{The Kakimizu complex of $11_{192}.$ The knot $11_{192}$ is a 2-bridge knot.}
  \label{}
\end{figure}

The list of Kakimizu complex of $2$-bridge knots of $11$ crossings are as follows:

\begin{itemize}
    \item $$K = 11_{13}, 11_{59}, 11_{65}, 11_{75}, 11_{77}, 11_{84}, 11_{85}, 11_{89}, 11_{90}, 11_{91}, 11_{93},  11_{110}, 11_{111},$$ $$ 11_{117}, 11_{120}, 11_{140}, 11_{144}, 11_{178}, 11_{183}, 11_{185}, 11_{188}, 11_{190}, 11_{193}, 11_{195}, 11_{204}, 11_{205}, $$ $$ 11_{207}, 11_{208}, 11_{211}, 11_{220}, 11_{224}, 11_{225}, 11_{230}, 11_{234}, 11_{242}, 11_{246}, 11_{247}, 11_{307}, 11_{309}, $$ $$ 11_{334}, 11_{339}, 11_{342}, 11_{355}, 11_{358}, 11_{364}.$$

    are spanned by a unique minimal genus Seifert surface up to isotopy.

The Kakimizu complex of the knots are a single vertex.

    \begin{tikzpicture}
    \draw[olive,thick,latex-latex] (0,0) 
    node[pos=0,mynode,fill=red,label=above:\textcolor{red}{$T$}]{};
  \end{tikzpicture}

  \item The list of $2$-bridge knots with non-trivial Kakimizu complex.

\begin{table}[H]
\begin{tabular}{|l|l|l|l|}
\hline
 Knots & 2-bridge index & Even continued fraction expansion & Kakimizu complex \\ \hline
$11_{95}$ & $33/73 = -40/73$ & $[-2,-6,-4,-2]$ &     \begin{tikzpicture}
    \draw[olive,thick,latex-latex] (0,0) -- (1,0)
    node[pos=0,mynode,fill=red,label=above:\textcolor{red}{$T_1$}]{}
    node[pos=1,mynode,fill=blue,text=blue,label=above:\textcolor{blue}{$T_2$}]{};
  \end{tikzpicture}
 \\ \hline
 
 $11_{98}$& $18/47$ & $ [4,-4,-2,2]$ &      \begin{tikzpicture}
    \draw[olive,thick,latex-latex] (0,0) -- (1,0)
    node[pos=0,mynode,fill=red,label=above:\textcolor{red}{$T_1$}]{}
    node[pos=1,mynode,fill=blue,text=blue,label=above:\textcolor{blue}{$T_2$}]{};
  \end{tikzpicture}
\\ \hline
 
$11_{119}$ & $64/109$ & $ [2,-4,-4,2]$, &      \begin{tikzpicture}
    \draw[olive,thick,latex-latex] (0,0) -- (1,0)
    node[pos=0,mynode,fill=red,label=above:\textcolor{red}{$T_1$}]{}
    node[pos=1,mynode,fill=blue,text=blue,label=above:\textcolor{blue}{$T_2$}]{};
  \end{tikzpicture}
\\ \hline
 $11_{145}$ & $22/83$ & $ [4,4,-2,2]$, &  \begin{tikzpicture}
    \draw[olive,thick,latex-latex] (0,0) -- (1,0)
    node[pos=0,mynode,fill=red,label=above:\textcolor{red}{$T_1$}]{}
    node[pos=1,mynode,fill=blue,text=blue,label=above:\textcolor{blue}{$T_2$}]{};
  \end{tikzpicture}
  \\ \hline
$11_{154}$ & $37/67 = -30/67$ & $ [-2,4,-4,-2]$ &  \begin{tikzpicture}
    \draw[olive,thick,latex-latex] (0,0) -- (1,0)
    node[pos=0,mynode,fill=red,label=above:\textcolor{red}{$T_1$}]{}
    node[pos=1,mynode,fill=blue,text=blue,label=above:\textcolor{blue}{$T_2$}]{};
  \end{tikzpicture}\\ \hline
 
$11_{166}$ & $45/59 = -14/59$ & $ [-4,4,-2,-2]$ &  \begin{tikzpicture}
    \draw[olive,thick,latex-latex] (0,0) -- (1,0)
    node[pos=0,mynode,fill=red,label=above:\textcolor{red}{$T_1$}]{}
    node[pos=1,mynode,fill=blue,text=blue,label=above:\textcolor{blue}{$T_2$}]{};
  \end{tikzpicture}\\ \hline

$11_{186}$ & $39/95 = -56/95 $ & $ [-2,-4,-2,-2,-4,-2]$ &   \begin{tikzpicture}
    \draw[olive,thick,latex-latex] (0,0) -- (1,0)
    node[pos=0,mynode,fill=red,label=above:\textcolor{red}{$T_1$}]{}
    node[pos=1,mynode,fill=blue,text=blue,label=above:\textcolor{blue}{$T_2$}]{};

    \draw[olive,thick,latex-latex] (1,0) -- (2,0)
    node[pos=0,mynode,fill=blue,label=above:\textcolor{blue}{$T_2$}]{}
    node[pos=1,mynode,fill=green,text=green,label=above:\textcolor{green}{$T_3$}]{};

    \draw[olive,thick,latex-latex] (2,0) -- (3,0)
    node[pos=0,mynode,fill=green,label=above:\textcolor{green}{$T_3$}]{}
    node[pos=1,mynode,fill=brown,text=brown,label=above:\textcolor{brown}{$T_4$}]{};
  \end{tikzpicture}
\\ \hline

$11_{191}$ & $19/83 = -64/83$ & $ [-2,-2,-2,-4,-4,-2]$, &  \begin{tikzpicture}
    \draw[olive,thick,latex-latex] (0,0) -- (1,0)
    node[pos=0,mynode,fill=red,label=above:\textcolor{red}{$T_1$}]{}
    node[pos=1,mynode,fill=blue,text=blue,label=above:\textcolor{blue}{$T_2$}]{};
  \end{tikzpicture}\\ \hline
$11_{192}$ & $45/59 = -14/59$ & $ [-4,4,-2,-2]$ &  \begin{tikzpicture}

    \fill[fill=black!20] (0,0)--(1,0)--(1,1);

    \fill[fill=red!20] (1,0)--(2,0)--(1,1);

    \draw[olive,thick,latex-latex] (1,1) -- (1,0)
    node[pos=0,mynode,fill=green,label=above:\textcolor{green}{$T_3$}]{}
    node[pos=1,mynode,fill=blue,text=blue,label=above:\textcolor{blue}{$T_2$}]{};

    \draw[olive,thick,latex-latex] (0,0) -- (1,0)
    node[pos=0,mynode,fill=red,label=above:\textcolor{red}{$T_1$}]{}
    node[pos=1,mynode,fill=blue,text=blue,label=above:\textcolor{blue}{$T_2$}]{};

  \draw[olive,thick,latex-latex] (0,0) -- (1,1)
    node[pos=0,mynode,fill=red,label=above:\textcolor{red}{$T_1$}]{}
    node[pos=1,mynode,fill=green,text=blue,label=above:\textcolor{green}{$T_3$}]{};

    \draw[olive,thick,latex-latex] (1,1) -- (1,0)
    node[pos=0,mynode,fill=green,label=above:\textcolor{green}{$T_3$}]{}
    node[pos=1,mynode,fill=blue,text=blue,label=above:\textcolor{blue}{$T_2$}]{};

    \draw[olive,thick,latex-latex] (2,0) -- (1,0)
    node[pos=0,mynode,fill=brown,label=above:\textcolor{brown}{$T_4$}]{}
    node[pos=1,mynode,fill=blue,text=blue,label=above:\textcolor{blue}{$T_2$}]{};

  \draw[olive,thick,latex-latex] (2,0) -- (1,1)
    node[pos=0,mynode,fill=brown,label=above:\textcolor{brown}{$T_4$}]{}
    node[pos=1,mynode,fill=green,text=blue,label=above:\textcolor{green}{$T_3$}]{};
  \end{tikzpicture}

  \\ \hline
\end{tabular}
\end{table}

\begin{table}[H]
\begin{tabular}{|l|l|l|l|}
\hline
 Knots & 2-bridge index & Continued fraction expansion & Kakimizu complex
\\ \hline
$11_{210}$ & $16/73$ & $ [4,-2,-4,-2]$, &  \begin{tikzpicture}
    \draw[olive,thick,latex-latex] (0,0) -- (1,0)
    node[pos=0,mynode,fill=red,label=above:\textcolor{red}{$T_1$}]{}
    node[pos=1,mynode,fill=blue,text=blue,label=above:\textcolor{blue}{$T_2$}]{};

    \draw[olive,thick,latex-latex] (1,0) -- (2,0)
    node[pos=0,mynode,fill=blue,label=above:\textcolor{blue}{$T_2$}]{}
    node[pos=1,mynode,fill=green,text=green,label=above:\textcolor{green}{$T_3$}]{};
  \end{tikzpicture}
\\ \hline
$11_{226}$ & $20/71$ & $ [4,2,-4,2]$, &  \begin{tikzpicture}
    \draw[olive,thick,latex-latex] (0,0) -- (1,0)
    node[pos=0,mynode,fill=red,label=above:\textcolor{red}{$T_1$}]{}
    node[pos=1,mynode,fill=blue,text=blue,label=above:\textcolor{blue}{$T_2$}]{};

    \draw[olive,thick,latex-latex] (1,0) -- (2,0)
    node[pos=0,mynode,fill=blue,label=above:\textcolor{blue}{$T_2$}]{}
    node[pos=1,mynode,fill=green,text=green,label=above:\textcolor{green}{$T_3$}]{};
  \end{tikzpicture}
\\ \hline 
$11_{229}$ & $55/71 = -16/71$ & $ [-4,2,-4,-2]$, &  \begin{tikzpicture}
    \draw[olive,thick,latex-latex] (0,0) -- (1,0)
    node[pos=0,mynode,fill=red,label=above:\textcolor{red}{$T_1$}]{}
    node[pos=1,mynode,fill=blue,text=blue,label=above:\textcolor{blue}{$T_2$}]{};

    \draw[olive,thick,latex-latex] (1,0) -- (2,0)
    node[pos=0,mynode,fill=blue,label=above:\textcolor{blue}{$T_2$}]{}
    node[pos=1,mynode,fill=green,text=green,label=above:\textcolor{green}{$T_3$}]{};
  \end{tikzpicture}
\\ \hline 

$11_{235}$ & $49/71 = -22/71$ & $ [-4,-2,-2,-2,-4,-2]$, &    \begin{tikzpicture}
    \draw[olive,thick,latex-latex] (0,0) -- (1,0)
    node[pos=0,mynode,fill=red,label=above:\textcolor{red}{$T_1$}]{}
    node[pos=1,mynode,fill=blue,text=blue,label=above:\textcolor{blue}{$T_2$}]{};

    \draw[olive,thick,latex-latex] (1,0) -- (2,0)
    node[pos=0,mynode,fill=blue,label=above:\textcolor{blue}{$T_2$}]{}
    node[pos=1,mynode,fill=green,text=green,label=above:\textcolor{green}{$T_3$}]{};

    \draw[olive,thick,latex-latex] (2,0) -- (3,0)
    node[pos=0,mynode,fill=green,label=above:\textcolor{green}{$T_3$}]{}
    node[pos=1,mynode,fill=brown,text=brown,label=above:\textcolor{brown}{$T_4$}]{};

    \draw[olive,thick,latex-latex] (3,0) -- (4,0)
    node[pos=0,mynode,fill=brown,label=above:\textcolor{brown}{$T_4$}]{}
    node[pos=1,mynode,fill=violet,text=violet,label=above:\textcolor{violet}{$T_5$}]{};
  \end{tikzpicture}
\\ \hline 

$11_{236}$ & $29/99 = -70/99$ & $ [-2,-2,-4,-2,-4,-2]$, &  \begin{tikzpicture}
    \draw[olive,thick,latex-latex] (0,0) -- (1,0)
    node[pos=0,mynode,fill=red,label=above:\textcolor{red}{$T_1$}]{}
    node[pos=1,mynode,fill=blue,text=blue,label=above:\textcolor{blue}{$T_2$}]{};

    \draw[olive,thick,latex-latex] (1,0) -- (2,0)
    node[pos=0,mynode,fill=blue,label=above:\textcolor{blue}{$T_2$}]{}
    node[pos=1,mynode,fill=green,text=green,label=above:\textcolor{green}{$T_3$}]{};
  \end{tikzpicture}
\\ \hline 

$11_{238}$ & $53/65 = -12/65$ & $ [-6,-2,-4,-2]$, &  \begin{tikzpicture}
    \draw[olive,thick,latex-latex] (0,0) -- (1,0)
    node[pos=0,mynode,fill=red,label=above:\textcolor{red}{$T_1$}]{}
    node[pos=1,mynode,fill=blue,text=blue,label=above:\textcolor{blue}{$T_2$}]{};

    \draw[olive,thick,latex-latex] (1,0) -- (2,0)
    node[pos=0,mynode,fill=blue,label=above:\textcolor{blue}{$T_2$}]{}
    node[pos=1,mynode,fill=green,text=green,label=above:\textcolor{green}{$T_3$}]{};
  \end{tikzpicture}
\\ \hline 
$11_{243}$ & $49/69 = -20/69$ & $ [-4,-2,-6,-2]$, &  \begin{tikzpicture}
    \draw[olive,thick,latex-latex] (0,0) -- (1,0)
    node[pos=0,mynode,fill=red,label=above:\textcolor{red}{$T_1$}]{}
    node[pos=1,mynode,fill=blue,text=blue,label=above:\textcolor{blue}{$T_2$}]{};

    \draw[olive,thick,latex-latex] (1,0) -- (2,0)
    node[pos=0,mynode,fill=blue,label=above:\textcolor{blue}{$T_2$}]{}
    node[pos=1,mynode,fill=green,text=green,label=above:\textcolor{green}{$T_3$}]{};
  \end{tikzpicture}
\\ \hline 
$11_{311}$ & $61/79 = -18/79$ & $ [-4,2,-2,-4]$, &    \begin{tikzpicture}
    \draw[olive,thick,latex-latex] (0,0) -- (1,0)
    node[pos=0,mynode,fill=red,label=above:\textcolor{red}{$T_1$}]{}
    node[pos=1,mynode,fill=blue,text=blue,label=above:\textcolor{blue}{$T_2$}]{};

    \draw[olive,thick,latex-latex] (1,0) -- (2,0)
    node[pos=0,mynode,fill=blue,label=above:\textcolor{blue}{$T_2$}]{}
    node[pos=1,mynode,fill=green,text=green,label=above:\textcolor{green}{$T_3$}]{};

    \draw[olive,thick,latex-latex] (2,0) -- (3,0)
    node[pos=0,mynode,fill=green,label=above:\textcolor{green}{$T_3$}]{}
    node[pos=1,mynode,fill=brown,text=brown,label=above:\textcolor{brown}{$T_4$}]{};
  \end{tikzpicture}

\\ \hline 
$11_{333}$ & $14/65$ & $ [4,-2,-2,4]$, &    \begin{tikzpicture}
    \draw[olive,thick,latex-latex] (0,0) -- (1,0)
    node[pos=0,mynode,fill=red,label=above:\textcolor{red}{$T_1$}]{}
    node[pos=1,mynode,fill=blue,text=blue,label=above:\textcolor{blue}{$T_2$}]{};

    \draw[olive,thick,latex-latex] (1,0) -- (2,0)
    node[pos=0,mynode,fill=blue,label=above:\textcolor{blue}{$T_2$}]{}
    node[pos=1,mynode,fill=green,text=green,label=above:\textcolor{green}{$T_3$}]{};

    \draw[olive,thick,latex-latex] (2,0) -- (3,0)
    node[pos=0,mynode,fill=green,label=above:\textcolor{green}{$T_3$}]{}
    node[pos=1,mynode,fill=brown,text=brown,label=above:\textcolor{brown}{$T_4$}]{};
  \end{tikzpicture}

\\ \hline 
$11_{335}$ & $17/75 = -58/75$ & $ [-2,-2,-2,-4,-2,-4]$, &  \begin{tikzpicture}
    \draw[olive,thick,latex-latex] (0,0) -- (1,0)
    node[pos=0,mynode,fill=red,label=above:\textcolor{red}{$T_1$}]{}
    node[pos=1,mynode,fill=blue,text=blue,label=above:\textcolor{blue}{$T_2$}]{};

    \draw[olive,thick,latex-latex] (1,0) -- (2,0)
    node[pos=0,mynode,fill=blue,label=above:\textcolor{blue}{$T_2$}]{}
    node[pos=1,mynode,fill=green,text=green,label=above:\textcolor{green}{$T_3$}]{};
  \end{tikzpicture}
\\ \hline 

$11_{336}$ & $11/59 = -48/59$ & $ [-2,-2,-2,-2,-4,-4]$, &  \begin{tikzpicture}
    \draw[olive,thick,latex-latex] (0,0) -- (1,0)
    node[pos=0,mynode,fill=red,label=above:\textcolor{red}{$T_1$}]{}
    node[pos=1,mynode,fill=blue,text=blue,label=above:\textcolor{blue}{$T_2$}]{};
  \end{tikzpicture}\\ \hline
$11_{337}$ & $63/89 = -16/89$ & $ [-6,-2,4,2]$, &  \begin{tikzpicture}
    \draw[olive,thick,latex-latex] (0,0) -- (1,0)
    node[pos=0,mynode,fill=red,label=above:\textcolor{red}{$T_1$}]{}
    node[pos=1,mynode,fill=blue,text=blue,label=above:\textcolor{blue}{$T_2$}]{};

    \draw[olive,thick,latex-latex] (1,0) -- (2,0)
    node[pos=0,mynode,fill=blue,label=above:\textcolor{blue}{$T_2$}]{}
    node[pos=1,mynode,fill=green,text=green,label=above:\textcolor{green}{$T_3$}]{};
  \end{tikzpicture}
\\ \hline 
$11_{343}$ & $27/31 = -4/31$ & $ [-8,-4]$, &  \begin{tikzpicture}
    \draw[olive,thick,latex-latex] (0,0) -- (1,0)
    node[pos=0,mynode,fill=red,label=above:\textcolor{red}{$T_1$}]{}
    node[pos=1,mynode,fill=blue,text=blue,label=above:\textcolor{blue}{$T_2$}]{};
  \end{tikzpicture}\\ \hline
$11_{356}$ & $55/79 = -24/79$ & $ [-4,-2,-2,-4,-2,-2]$, &    \begin{tikzpicture}
    \draw[olive,thick,latex-latex] (0,0) -- (1,0)
    node[pos=0,mynode,fill=red,label=above: \textcolor{red}{$T_1$}]{}
    node[pos=1,mynode,fill=blue,text=blue,label=above:\textcolor{blue}{$T_2$}]{};

    \draw[olive,thick,latex-latex] (1,0) -- (2,0)
    node[pos=0,mynode,fill=blue,label=above:\textcolor{blue}{$T_2$}]{}
    node[pos=1,mynode,fill=green,text=green,label=above:\textcolor{green}{$T_3$}]{};

    \draw[olive,thick,latex-latex] (2,0) -- (3,0)
    node[pos=0,mynode,fill=green,label=above:\textcolor{green}{$T_3$}]{}
    node[pos=1,mynode,fill=brown,text=brown,label=above:\textcolor{brown}{$T_4$}]{};
  \end{tikzpicture}

\\ \hline 

$11_{357}$ & $27/91 = -64/91$ & $ [-2,-2,-4,-4,-2,-2]$, &  \begin{tikzpicture}
    \draw[olive,thick,latex-latex] (0,0) -- (1,0)
    node[pos=0,mynode,fill=red,label=above:\textcolor{red}{$T_1$}]{}
    node[pos=1,mynode,fill=blue,text=blue,label=above:\textcolor{blue}{$T_2$}]{};
  \end{tikzpicture}\\ \hline

$11_{359}$ & $43/53 = -10/79$ & $ [-6,-2,-2,-4]$, &    \begin{tikzpicture}
    \draw[olive,thick,latex-latex] (0,0) -- (1,0)
    node[pos=0,mynode,fill=red,label=above: \textcolor{red}{$T_1$}]{}
    node[pos=1,mynode,fill=blue,text=blue,label=above:\textcolor{blue}{$T_2$}]{};

    \draw[olive,thick,latex-latex] (1,0) -- (2,0)
    node[pos=0,mynode,fill=blue,label=above:\textcolor{blue}{$T_2$}]{}
    node[pos=1,mynode,fill=green,text=green,label=above:\textcolor{green}{$T_3$}]{};

    \draw[olive,thick,latex-latex] (2,0) -- (3,0)
    node[pos=0,mynode,fill=green,label=above:\textcolor{green}{$T_3$}]{}
    node[pos=1,mynode,fill=brown,text=brown,label=above:\textcolor{brown}{$T_4$}]{};
  \end{tikzpicture}

\\ \hline

$11_{360}$ & $47/57 = -10/59$ & $ [-6,-4,-2,-2]$, &  \begin{tikzpicture}
    \draw[olive,thick,latex-latex] (0,0) -- (1,0)
    node[pos=0,mynode,fill=red,label=above:\textcolor{red}{$T_1$}]{}
    node[pos=1,mynode,fill=blue,text=blue,label=above:\textcolor{blue}{$T_2$}]{};
  \end{tikzpicture}\\ \hline

$11_{363}$ & $29/35 = -6/35$ & $ [-6,-6]$, &  \begin{tikzpicture}
    \draw[olive,thick,latex-latex] (0,0) -- (1,0)
    node[pos=0,mynode,fill=red,label=above:\textcolor{red}{$T_1$}]{}
    node[pos=1,mynode,fill=blue,text=blue,label=above:\textcolor{blue}{$T_2$}]{};
  \end{tikzpicture}\\ \hline

$11_{365}$ & $35/51 = -16/51$ & $ [-4,-2,-2,-2,-2,-4]$, &  \begin{tikzpicture}
    \draw[olive,thick,latex-latex] (0,0) -- (1,0)
    node[pos=0,mynode,fill=red,label=above:\textcolor{red}{$T_1$}]{}
    node[pos=1,mynode,fill=blue,text=blue,label=above:\textcolor{blue}{$T_2$}]{};

    \draw[olive,thick,latex-latex] (1,0) -- (2,0)
    node[pos=0,mynode,fill=blue,label=above:\textcolor{blue}{$T_2$}]{}
    node[pos=1,mynode,fill=green,text=green,label=above:\textcolor{green}{$T_3$}]{};

    \draw[olive,thick,latex-latex] (2,0) -- (3,0)
    node[pos=0,mynode,fill=green,label=above:\textcolor{green}{$T_3$}]{}
    node[pos=1,mynode,fill=brown,text=brown,label=above:\textcolor{brown}{$T_4$}]{};

    \draw[olive,thick,latex-latex] (3,0) -- (4,0)
    node[pos=0,mynode,fill=brown,label=above:\textcolor{brown}{$T_4$}]{}
    node[pos=1,mynode,fill=violet,text=violet,label=above:\textcolor{violet}{$T_5$}]{};

    \draw[olive,thick,latex-latex] (4,0) -- (5,0)
    node[pos=0,mynode,fill=violet,label=above:\textcolor{violet}{$T_5$}]{}
    node[pos=1,mynode,fill=orange,text=orange,label=above:\textcolor{orange}{$T_6$}]{};
  \end{tikzpicture}
\\ \hline
\end{tabular}
\end{table}

\end{itemize}

\clearpage

\subsection{General  Algorithm}

Let $D$ be a reduced oriented alternating diagram of an 11 crossing  prime, alternating knot $K.$ 

\begin{itemize}

\item Firstly we identify (from the table \cite{knotinfo}) if the given $11$ crossing alternating link is a $2$-bridge knot. We compute the Kakimizu complex using Hatcher and Thurston's \cite{Hatcher1985} result to find $MS(L).$ We use the generalised result to find the $1$-skeleton of the Kakimizu complex of $2$-bridge knots.

    \item If it is not a $2$-bridge knot, we can apply the Seifert's algorithm on $D.$ We get a minimal genus Seifert surface $S$ with $\partial S = K.$

    We know that the Kakimizu complex is a flag complex. Hence finding the graph containing the $0$ and $1$ skeleton of the Kakimizu complex suffices to find the whole complex. 

    Secondly, the Kakimizu complex is connected, implies that starting from a minimal genus Seifert surface $S$ on the knot $K,$ we need to obtain all vertices adjacent to $S$ and repeat the process for every new vertex found.

  Let L be a non-split prime alternating link with $n > 0$ crossings.
    It has been shown by Hass,Thompson and Tsvietkova that for each fixed $g$ the number of genus $g-$ Seifert surfaces for $L$ is bounded by an explicitly given polynomial in $n.$ 
    
     Hence we know that the process in the previous part ends, in other words the Kakimizu complex for a link $L$ is finite.  

\item Any alternating link $K$ with a diagram $D$ can be written as a $*-$ product of special alternating links $K = *(L_i)_i$. When we apply Seifert's algorithm on an oriented alternating reduced diagram $D$ of $K$ then we get the surface $S$ as a Murasugi sums of spanning surfaces $S_i$ on special alternating links $L_i$.

\item If $K$ is a special alternating link, then we use Jessica Banks' result to find the Kakimizu complex of $K.$

\item Identify all the fibred pieces (if any) plumbed on a surface $S_1$ to obtain $S.$ 

\item If $S_1$ is non-fibred and has a unique Seifert surface, and is Murasugi summed with fibred surfaces, then the Kakimizu complex of $K$ is a single vertex.

\end{itemize}

\newpage

Let $L$ be a non-split, prime link and let $S$ be a minimal genus Seifert surface on $L.$ $IS(L,S)$ is the set of all isotopy classes of surfaces which can be made disjoint from $S$ in the link complement. (With respect to Kakimizu complex, surfaces(upto isotopy) in $IS(L,S)$ represents vertices in the Kakimizu complex which share an edge with $v_{[S]}$).    
    
    \vspace{.5cm}
    
    Boileau and Gabai showed that,
    
    \vspace{.5cm}
    
    \emph{Let $L$ be a non-split oriented link and $S$ a connected minimal genus Seifert surface for $L$. Suppose that $S$ is a Murasugi sum of $S_1$ and $S_2$, where $S_i$ is a spanning surface for an oriented link $L_i$ $(i=1,2)$.
Suppose further that $L_2$ is a fibred link with fibre $S_2$. Then $L_1$ is non-split, and
$S_1$ is connected and minimal genus. Moreover there is a bijection}

$$\phi:IS(L,S)\to IS(L_1,S_1).$$

 Let 
    
    $$S = S_1 \cup_{D_1} F_1 \cup_{D_2} F_2 \cup \dots \cup_{D_n}F_n $$  where 
    $F_i$ are fibred surfaces on $\partial F_{i}$ and each $S_1\cup_{D_i}F_i$ is a plumbing with $D_i$ being the plumbing disks in the surface $S_i.$ Note that all the plumbing disks $D_i$ are disjoint. Let $\partial S_1 = L_1$. Then $L_1$ is non-split and $S_1$ is connected.
    
    Assume $$IS(L_1) = [S_1].$$
    
    Then $$IS(L) = [S].$$

\begin{proof}

$IS(L_1,S_1) = \Phi.$

Consider $$S_1^1 = S_1\cup_{D_1}F_1.$$ We know that $F_1$ is a fibred surface for $L_1$. 
From the theorem by Boileau and Gabai, we have 

$$IS(\partial S_1^1, S_1^1) = \Phi.$$

Define $$S_1^k = S_1 \cup_{D_1} F_1 \cup_{D_2} F_2 \cup \dots \cup_{D_k}F_k $$ for every $1\leq k\leq n$

Repeating the same process as above we have, 

$$IS(\partial S_1^k, S_1^k) = \Phi.$$

For $k=n$ we have $$IS(L,S) = \Phi \implies IS(L) = [S].$$

\end{proof}

The list of the Kakimizu complexes of  $11$ crossing alternating knots $K$, that are Murasugi sums of links $L_1$ with a unique incompressible surface and $L_2$ that are fibred:

\begin{itemize}
    \item The following knots $K$ spans an unique minimal genus Seifert surface up to isotopy.
    
    $$K = 11_{1},11_{2},11_{4},11_{6},11_{8},11_{10},11_{11},11_{12},11_{16},11_{18},11_{20},11_{21},11_{23},11_{27},11_{29},11_{30},11_{31},$$
    $$11_{32},11_{36},11_{37},11_{38},11_{39},11_{41},11_{46},11_{48},11_{49},11_{50},11_{52},11_{54},11_{56},11_{58},11_{60},11_{63},11_{64},$$
    $$11_{67},11_{69},11_{70},11_{78},11_{87},11_{97},11_{100},11_{101},11_{102},11_{104},11_{105},11_{107},11_{114},11_{115},11_{116},$$
    $$11_{118},11_{122},11_{130},11_{132},11_{133},11_{134},11_{135},11_{136},11_{137},11_{138},11_{141},11_{143},11_{148},11_{149},11_{150},$$
    $$11_{152},11_{153},11_{155},11_{161},11_{165},11_{167},11_{168},11_{169},11_{172},11_{173},11_{181},11_{187},11_{197},11_{198},11_{199},$$
    $$11_{202},11_{212},11_{213},11_{214},11_{219},11_{222},11_{249},11_{256},11_{258},11_{260},11_{262},11_{265},11_{270},11_{271},11_{272},$$
    $$11_{273},11_{275},11_{276},11_{278},11_{279},11_{283},11_{285},11_{290},11_{294},11_{295},11_{296},11_{297},11_{303},11_{304},11_{312},$$
    $$11_{313},11_{317},11_{321},11_{322},11_{323},11_{324},11_{327},11_{328},11_{331},11_{344},11_{345},11_{347},11_{349},11_{352}.$$

The Kakimizu complex of the knot $K$ is: 

    \begin{tikzpicture}
    \draw[olive,thick,latex-latex] (0,0) 
    node[pos=0,mynode,fill=red,label=above:\textcolor{red}{$T$}]{};
  \end{tikzpicture}

\item Let $K$ be an $11$ crossing alternating knot. Let $D$ be a reduced, prime, alternating, oriented diagram of $L.$ Let $S$ be a Seifert surface obtained by applying Seifert's algorithm on $D.$

Let $S = S_1\cup_{D_0}S_2$ such that $S_2$ be a fibred surface on a fibred link $L_2.$ We have, 

$IS(L,S)\cong IS(L_1,S_1).$ Let us call this bijective map $$f:IS(L,S)\to IS(L_1,S_1).$$

Assume that for every $T\in IS(L_1,S_1),$ $f^{-1} (T) = T\cup_{E}S_2.$ Assume that this condition holds for every representative surface of a generic vertex $S$ on $L.$ Then $IS(L)\cong IS(L_1)$ and the Kakimizu complexes are identical.

\begin{proof}
Let $D$ be a reduced, alternating, oriented diagram of the knot $K.$ Let $S$ be the Seifert surface obtained by applying Seifert's algorithm on $D.$

Let $v$ be a generic vertex in $MS(L) = IS(L).$ Since the Kakimizu complex of a prime, non-split link is connected, there exists a path, $[[S],v_1,v_2,\dots,v_k = v]$ in the Kakimizu complex.

Since $S = S_1\cup_{D_0}S_2$ with $S_2$ fibred, we have 

$$f_1:IS(L,S) \to IS(L_1,S_1); \text{ with } S = S_1\cup S_2.$$ For every $v\in IS(L_1,S_1)$, there exists a representative surface $T$ such that $v = [T]$ and $f_1^{-1}([T]) = [T\cup S_2].$

$v_1\in IS(L,S),$ $f_1(v_1)\in IS(L_1,S_1).$ There exists a representative surface $T_1$ on $L_1$ such that $$f_1(v_1) = [T_1] \text{ and } v_1 = [T_1\cup S_2].$$

Define a map $f:IS(L) \to IS(L_1)$ with $f([S]) = [S_1],$ and $f(v_1) = [T_1].$

$$f_2:IS(L,T_1\cup S_2) \to IS(L_1,T_1).$$

$v_2\in IS(L,T_1\cup S_2),$ $f_2(v_2)\in IS(L_1,S_1).$ There exists a representative surface $T_2$ on $L_1$ such that $$f_2(v_2) = [T_2] \text{ and } v_2 = [T_2\cup S_2].$$

$f(v_2) = [T_2].$ 

Continuing the process, we obtain $f(v).$

Hence we define a one-one correspondence between $IS(L)$ and $IS(L_1)$ and so the Kakimizu complexes of $L$ are identical to $L_1.$

\end{proof}

The list of the Kakimizu complexes of  $11$ crossing alternating knots $K$, that are Murasugi sums of links $L_1$ with minimal genus Seifert surfaces and $L_2$ that is fibred with fibre $S_2$ with the condition that every surface on $K$ is a Murasugi sum of a surface on $L_1$ and $S_2.$

$K = 11_{61}.$

$IS(K) = \{[S_1],[S_2]\}$ and the Kakimizu complex is:         \begin{tikzpicture}
    \draw[olive,thick,latex-latex] (0,0) -- (1,0)
    node[pos=0,mynode,fill=red,label=above:\textcolor{red}{$[S]$}]{}
    node[pos=1,mynode,fill=blue,text=blue,label=above:\textcolor{blue}{$[S^c]$}]{};
  \end{tikzpicture}

\end{itemize}

\clearpage

\subsection*{Kakimizu complex of plumbings of links with unique spanning surfaces}

The following section is adopted from \cite{10.32917/hmj/1150922486} and proved by Kakimizu, in 1992. 

Let $(M,\gamma)$ be a sutured manifold. A marking $A$ is a properly embedded arc in $R(\gamma).$ A sutured manifold with a prescribed marking $(M,\gamma,A)$ is called a \emph{marked sutured manifold}. If there is a product disk $D$ in $(M,\gamma)$ with $A$ as an edge of $\partial D,$ then the sutured manifold $(M,\gamma)$ with the opposite edge $B$ is also a marked sutured manifold $(M,\gamma,B).$

\emph{Let $(M,\gamma,A)$ be a marked sutured manifold. Suppose that $M$ is irreducible and each component of $R(\gamma)$ is incompressible . If there is a product disk with $A$ as an edge, then the ambient isotopy types of product disks are unique.} 

For our purpose, if $L$ is non-split then $E(L)$ and complementary sutured manifold, $(M,\gamma)$ for a Seifert surface $S$ are irreducible 3-manifolds. If $S$ is incompressible then $R(\gamma)$ for the complementary sutured manifold $(M,\gamma)$ is incompressible. Let $S$ be a plumbing of $S_1$ and $S_2.$ Consider the complementary sutured manifolds $(M_1,\gamma_1)$ and $(M_2,\gamma_2)$ for $S_1$ and $S_2$ respectively. Let $D_0$ be the plumbing disk. $D_0\subset S_1$ is an embedded disk. Consider the core curve of the disk $D_0$ in $S_1.$ Since $M = E(L) - N(S),$ push out the core curve on the side where $S_2$ is attached. The push off of the core curve is denoted by the marking $A_1\subset R(\gamma_1).$ Hence we get a marked sutured manifold $(M_1,\gamma_1,A_1).$ Similarly we obtain a marked sutured manifold for $S_2$ denoted by $(M_2,\gamma_2,A_2).$  

\begin{theorem}
    Let $L$ be a non-split, prime, alternating link and let $D$ be a reduced alternating, oriented diagram $D$ of $L.$ Let $S$ be the Seifert surface obtained by applying Seifert's algorithm on the diagram $D.$ Let $S$ be a plumbing of $S_1$ and $S_2$ which are unique minimal genus Seifert surfaces for links $L_1$ amd $L_2$ respectively. Assume that $S_1$ and $S_2$ are not fibred. Let $(M_i,\gamma_i,A_i)$ and $(M_i,\gamma_i,A_i')$ $(i = 1,2)$ be the marked sutured manifolds for $S = S_1\cup_{D_0}S_2$ and $S^{c} = S_1\cup_{D_0^c}S_2,$ the dual of $S$ respectively. Then 

    \begin{itemize}
        \item $MS(L) = \{[S],[S^{c}]\}$ and the Kakimizu complex is 

        \begin{tikzpicture}
    \draw[olive,thick,latex-latex] (0,0) -- (1,0)
    node[pos=0,mynode,fill=red,label=above:\textcolor{red}{$[S]$}]{}
    node[pos=1,mynode,fill=blue,text=blue,label=above:\textcolor{blue}{$[S^c]$}]{};
  \end{tikzpicture}

       provided it satisfies the condition:

       \begin{itemize}
           \item there is no product disk with $A_1$ or $A_1'$ in $M_1.$

           \item  there is no product disk with $A_2$ or $A_2'$ in $M_2.$
       \end{itemize}

       \item Let us assume that there is a product disk in $(M_1,\gamma_1,A_1)$ with $A_1$ as an edge. Let $B_1$ be the opposite mark of the product disk. Let $S = S_1\cup_{D_0} S_2.$ Suppose $T = S_1\cup_{E_0}T_2$ is the plumbing with respect to the marking $B_1.$ Note that $S=T.$

       $MS(L) = \{[S],[S^{c}],[T^{c}]\}$ and the Kakimizu complex is 

      \begin{tikzpicture}
    \draw[olive,thick,latex-latex] (0,0) -- (1,0)
    node[pos=0,mynode,fill=red,label=above:\textcolor{red}{$[S^c]$}]{}
    node[pos=1,mynode,fill=blue,text=blue,label=above:\textcolor{blue}{$[S]$}]{};

    \draw[olive,thick,latex-latex] (1,0) -- (2,0)
    node[pos=0,mynode,fill=blue,label=above:\textcolor{blue}{$[S]$}]{}
    node[pos=1,mynode,fill=green,text=green,label=above:\textcolor{green}{$[T^c]$}]{};
  \end{tikzpicture}
  
  provided it satisfies the condition:

  \begin{itemize}
      \item there is no product disk with $A_1'$ or $B_1'$ in $M_1$

      \item there is no product disk with $A_2$ or $A_2'$ in $M_2.$
  \end{itemize}
    \end{itemize}
\end{theorem}

The list of the Kakimizu complexes of  $11$ crossing alternating knots $K$, that spans a surface $S$ that are plumbings of links with unique incompressible surfaces :

\begin{itemize}
    \item $K = 11_{45}.$ Let $S = S_1\cup_{D_0} S_2$ on links $L_1$ and $L_2.$

    The link $L_1$ is a special alternating link with a unique minimal genus Seifert surface $S_1$. $S_2$ is the $4$-half-twisted band, a unique minimal genus Seifert surface. There are no product disks with the markings as an edge. Hence the Kakimizu complex of $K$ is:

        \begin{tikzpicture}
    \draw[olive,thick,latex-latex] (0,0) -- (1,0)
    node[pos=0,mynode,fill=red,label=above:\textcolor{red}{$[S]$}]{}
    node[pos=1,mynode,fill=blue,text=blue,label=above:\textcolor{blue}{$[S^c]$}]{};
  \end{tikzpicture}

  \item $K = 11_{280}.$

  $S = S_1\cup_{D_0}S_2$ with $L_1 = (2,2,2)$ and $L_2 = (3,1,1)$ Pretzel knots. Both $S_1,S_2$ are unique minimal genus spanning surfaces for $L_1,L_2$ respectively.  They Kakimizu complex is:

       \begin{tikzpicture}
    \draw[olive,thick,latex-latex] (0,0) -- (1,0)
    node[pos=0,mynode,fill=red,label=above:\textcolor{red}{$[S]$}]{}
    node[pos=1,mynode,fill=blue,text=blue,label=above:\textcolor{blue}{$[S^c]$}]{};
  \end{tikzpicture}

  \item $K = 11_{325}.$ $S = S_1\cup_{D_0}S_2$ with $S_1$, a 4-half twisted band and $S_2$ is a unique spanning surface on a special alternating link $L_2.$ There is a product disk with the marking as an edge. The Kakimizu complex is:

      \begin{tikzpicture}
    \draw[olive,thick,latex-latex] (0,0) -- (1,0)
    node[pos=0,mynode,fill=red,label=above:\textcolor{red}{$T_1$}]{}
    node[pos=1,mynode,fill=blue,text=blue,label=above:\textcolor{blue}{$T_2$}]{};

    \draw[olive,thick,latex-latex] (1,0) -- (2,0)
    node[pos=0,mynode,fill=blue,label=above:\textcolor{blue}{$T_2$}]{}
    node[pos=1,mynode,fill=green,text=green,label=above:\textcolor{green}{$T_3$}]{};
  \end{tikzpicture}

\end{itemize}
\clearpage

Example:

\begin{figure}[H]
  \includegraphics[width=15cm]{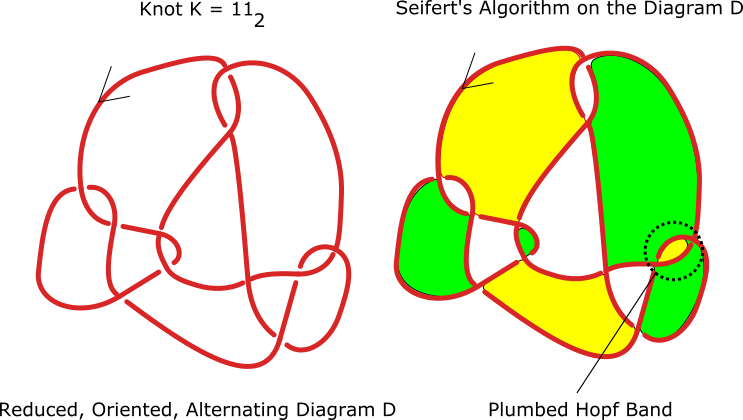}
  \caption{Knot $K$ = $11_2$ with a Seifert Surface on $K$}
  \label{fig:boat1}
\end{figure}

\begin{figure}[H]
  \includegraphics[width=15cm]{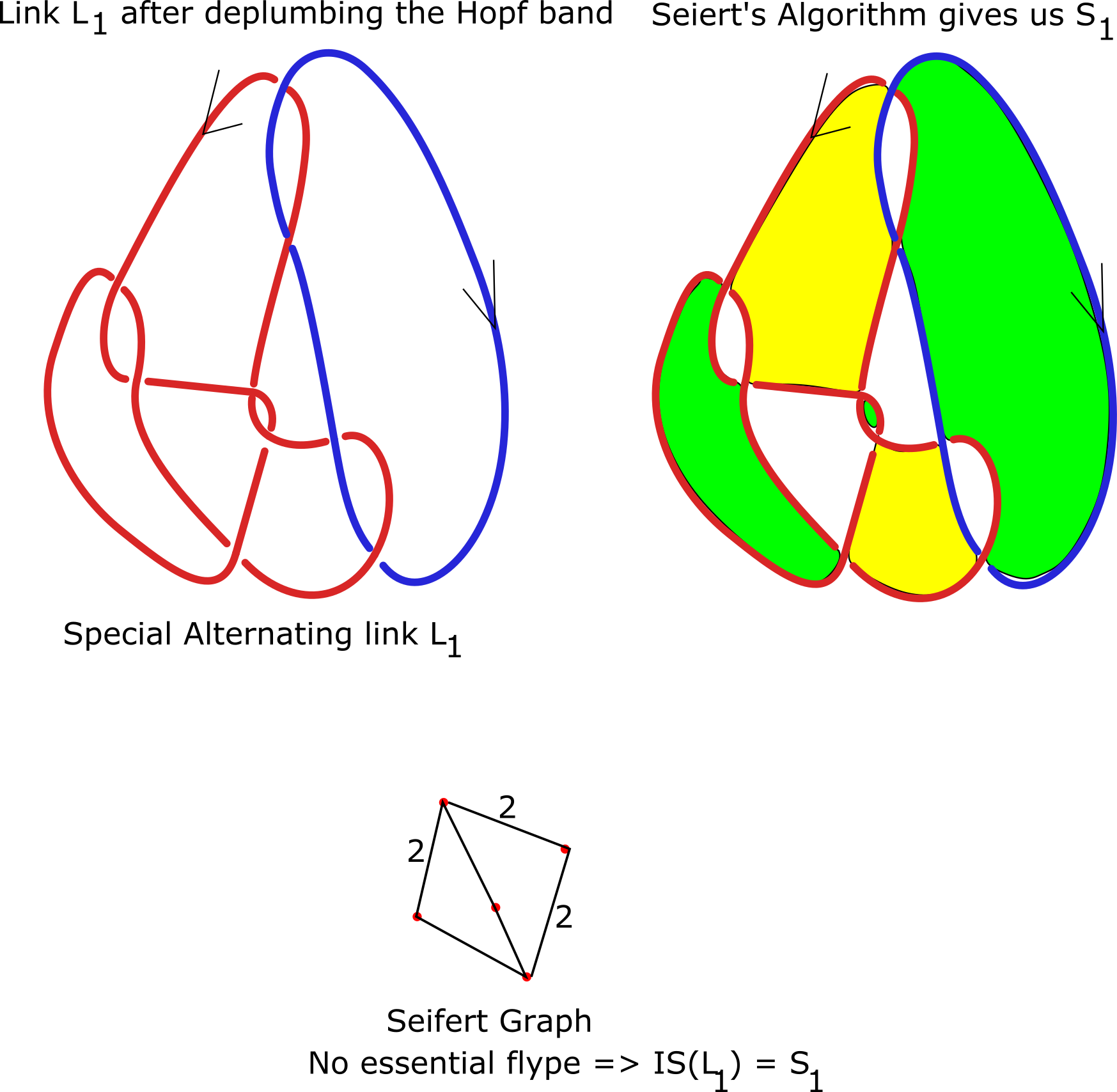}
  \caption{Deplumb the Hopf band to obtain a special alternating link $L_1$. There is a unique Seifert surface on $L_1.$}
  \label{fig:boat2}
\end{figure}

\newpage
\section*{Special case: $11_{103}$}

Let $D$ be a reduced oriented alternating diagram of an 
$11$ crossing prime alternating knot $K.$ We consider $K = 11_{103}.$

\begin{itemize}
    \item Apply Seifert's algorithm on the diagram $D$ of $K.$ We obtain a surface $T_1.$ 
    
    $$T_1 = S_1(7_4) \cup_{D_1}H_1\cup_{D_2}H_2$$
    
    where $S_1(7_4), S_2(7_4)$ are the two distinct surfaces (up to isotopy) as shown in the diagram and $H_i$ are Hopf bands plumbed onto $S_1(7_4).$
    
  \vspace{.5cm}
    
    Boileau and Gabai \cite{Gabai1986} showed that,
    
    \vspace{.5cm}
    
    \emph{Let $L$ be a non-split oriented link and $S$ a connected minimal genus Seifert surface for $L$. Suppose that $S$ is a Murasugi sum of $S_1$ and $S_2$, where $S_i$ is a spanning surface for an oriented link $L_i$ $(i=1,2)$.
Suppose further that $L_2$ is a fibred link with fibre $S_2$. Then $L_1$ is non-split, and
$S_1$ is connected and minimal genus. Moreover there is a bijection}

$$\phi:IS(L,S)\to IS(L_1,S_1)$$

Hence we have:

$$IS(11_{103},T_1) \xrightarrow{\phi_1}IS(L_0,S_1(7_4)\cup_{D_1}H_1)\xrightarrow{\phi_2}IS(7_4, S_1(7_4)) = [S_2(7_4)].$$
    
where $\phi_1$ and $\phi_2$ are isomorphisms. Hence we call $$IS(11_{103},T_1) = [T_2] = \phi_1^{-1}\circ \phi_2^{-1}([S_2(7_4)]).$$    

This gives us the following Lemma.

\begin{lemma}

$$IS(11_{103},T_1) = [T_2]$$ \emph{where $T_2$ is a non-plumbed disjoint surface in $E(L),$ not isotopic to $T_1.$ }

\end{lemma}

$$T_1 = S_1(7_4) \cup_{D_1} H_1 \cup_{D_2} H_2.$$

Since $H_1$ and $H_2$ are Hopf Bands, hence there are product disks passing through the Hopf bands. 

Product decomposition on $H_1$ and $H_2$ leads us to the complementary sutured manifold of $S_1(7_4).$ $IS(7_4,S_1(7_4)) = [S_2(7_4)].$ Inverting the operation of product decomposition on $H_1$ and $H_2$ yields us
$T_2$ in the complementary sutured manifold of $T_1.$

\end{itemize}

\begin{figure}[H]
  \includegraphics[width=16cm]{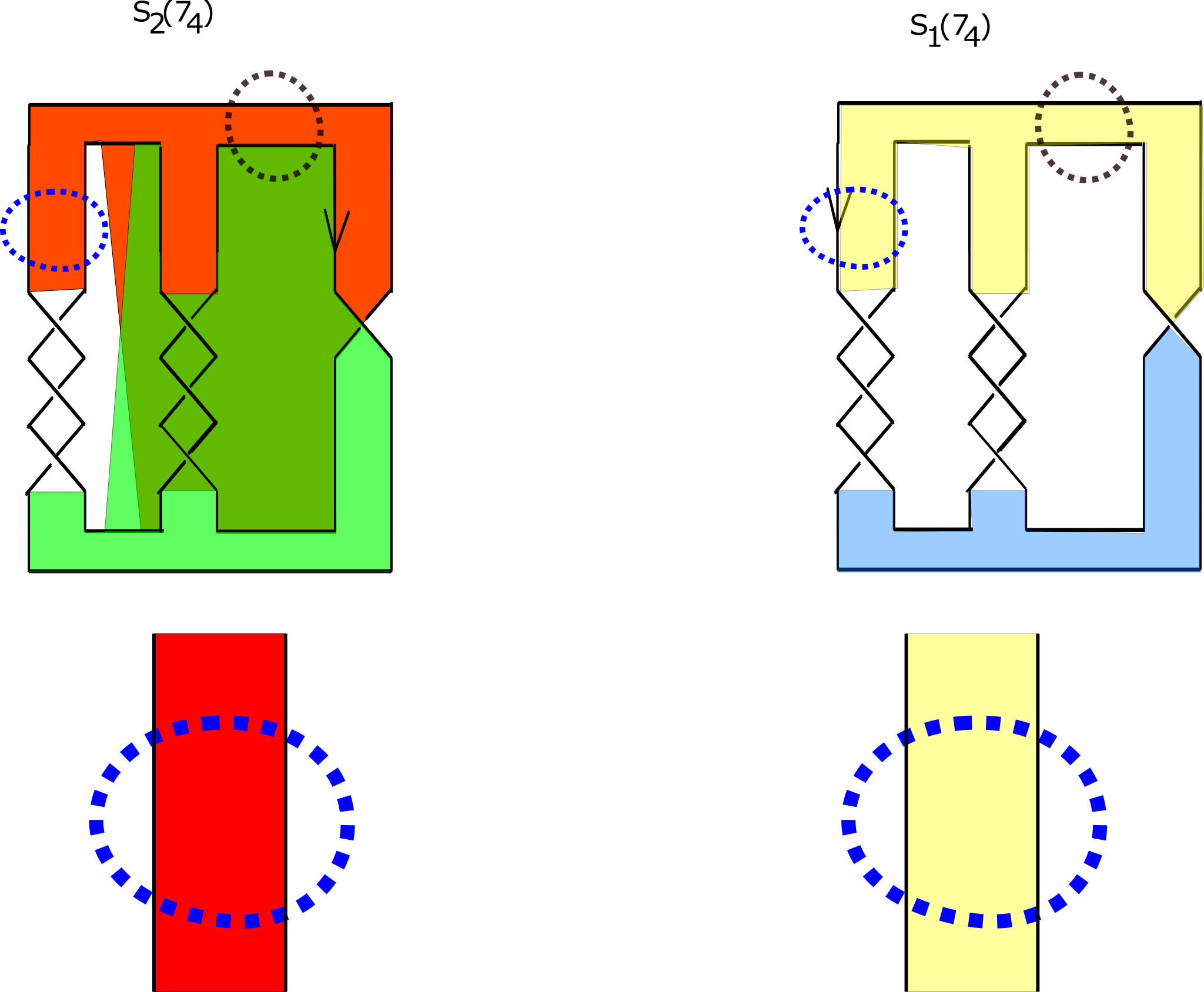}
  \caption{}
  \label{fig:boat1}
\end{figure}

\begin{figure}[H]
  \includegraphics[width=16cm]{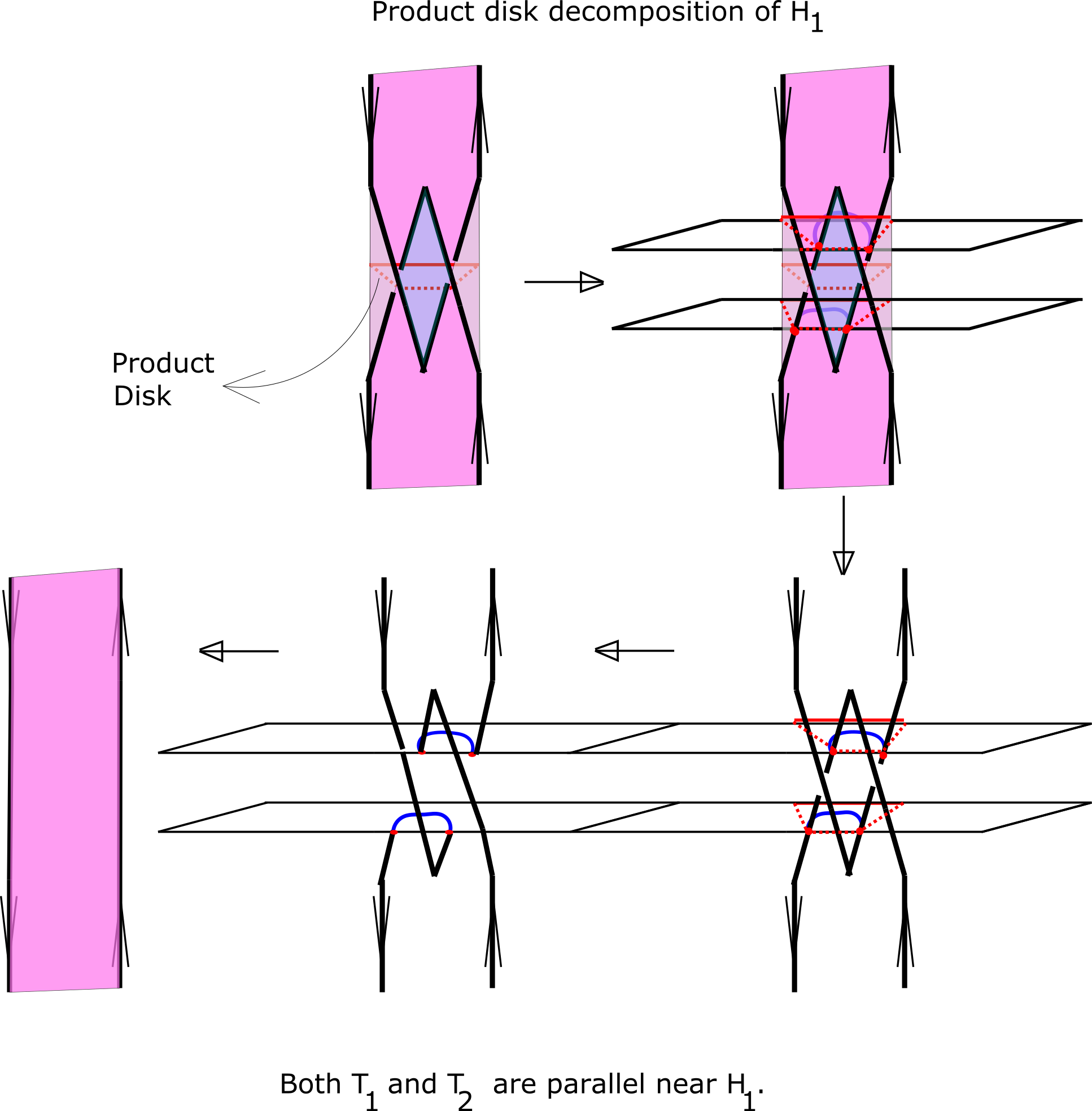}
  \caption{}
  \label{fig:boat1}
\end{figure}

\begin{figure}[H]
 \includegraphics[width=16cm]{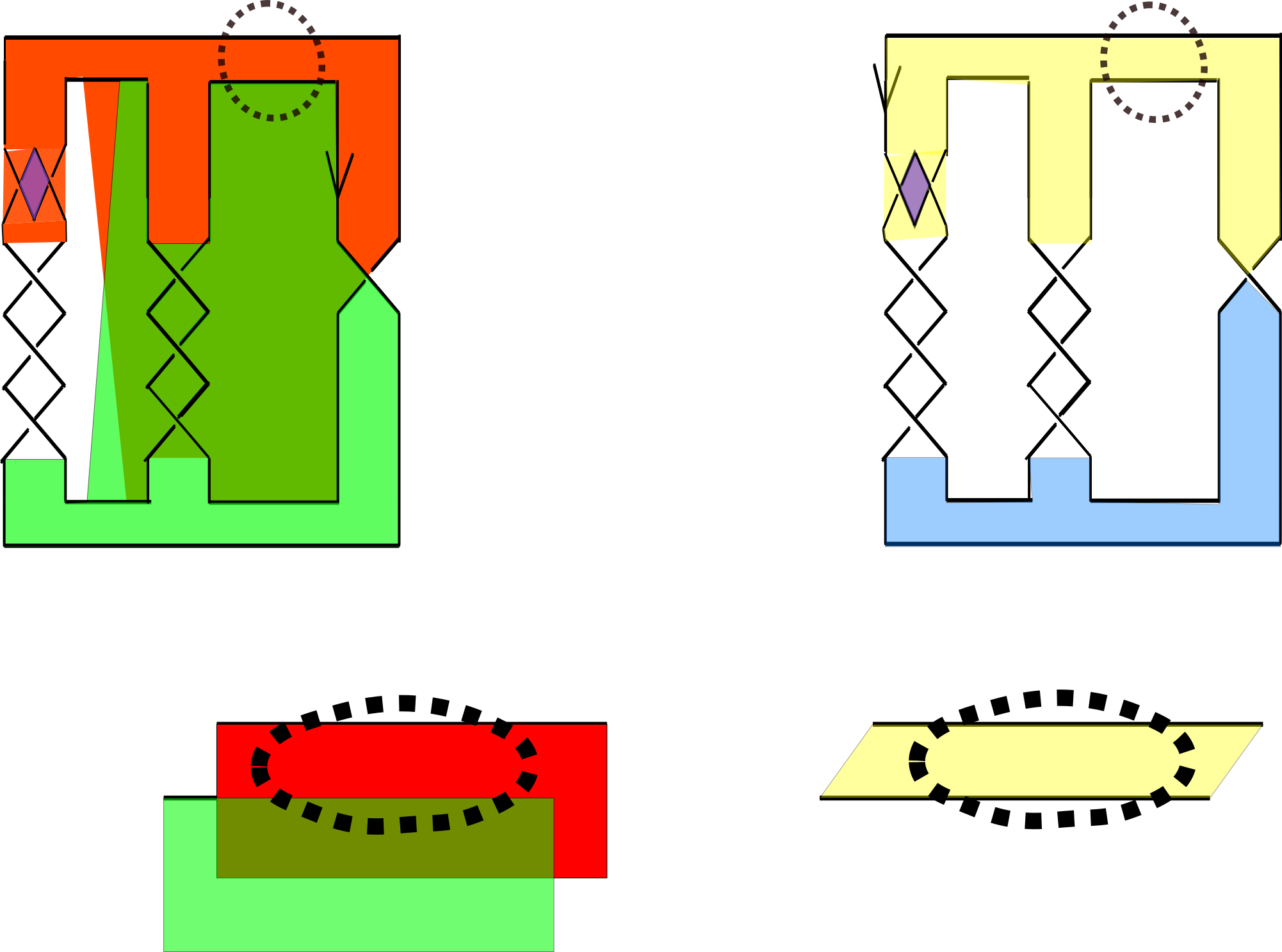}
  \caption{}
  \label{fig:boat1}
\end{figure}

\newpage

\begin{figure}[H]
  \includegraphics[width=9cm]{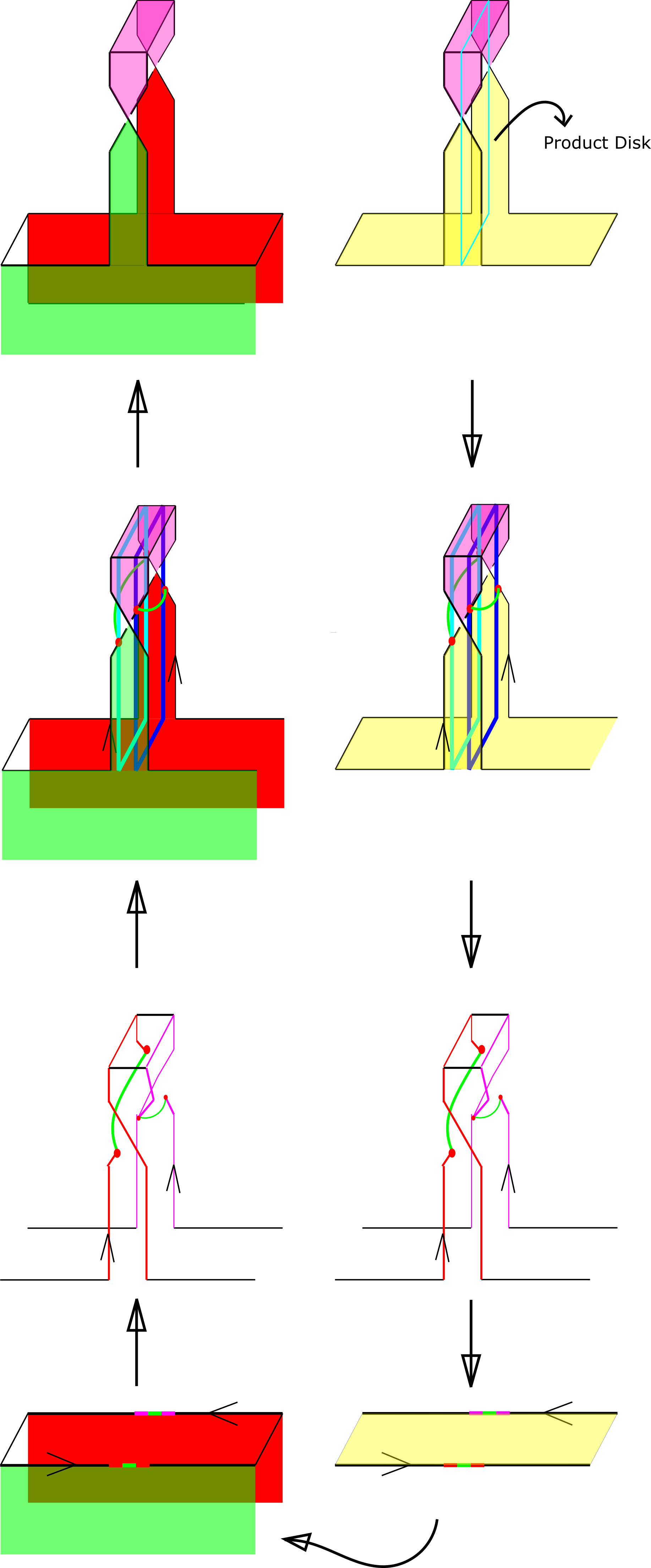}
  \caption{}
  \label{fig:boat1}
\end{figure}

\begin{figure}[H]
  \includegraphics[width=16cm]{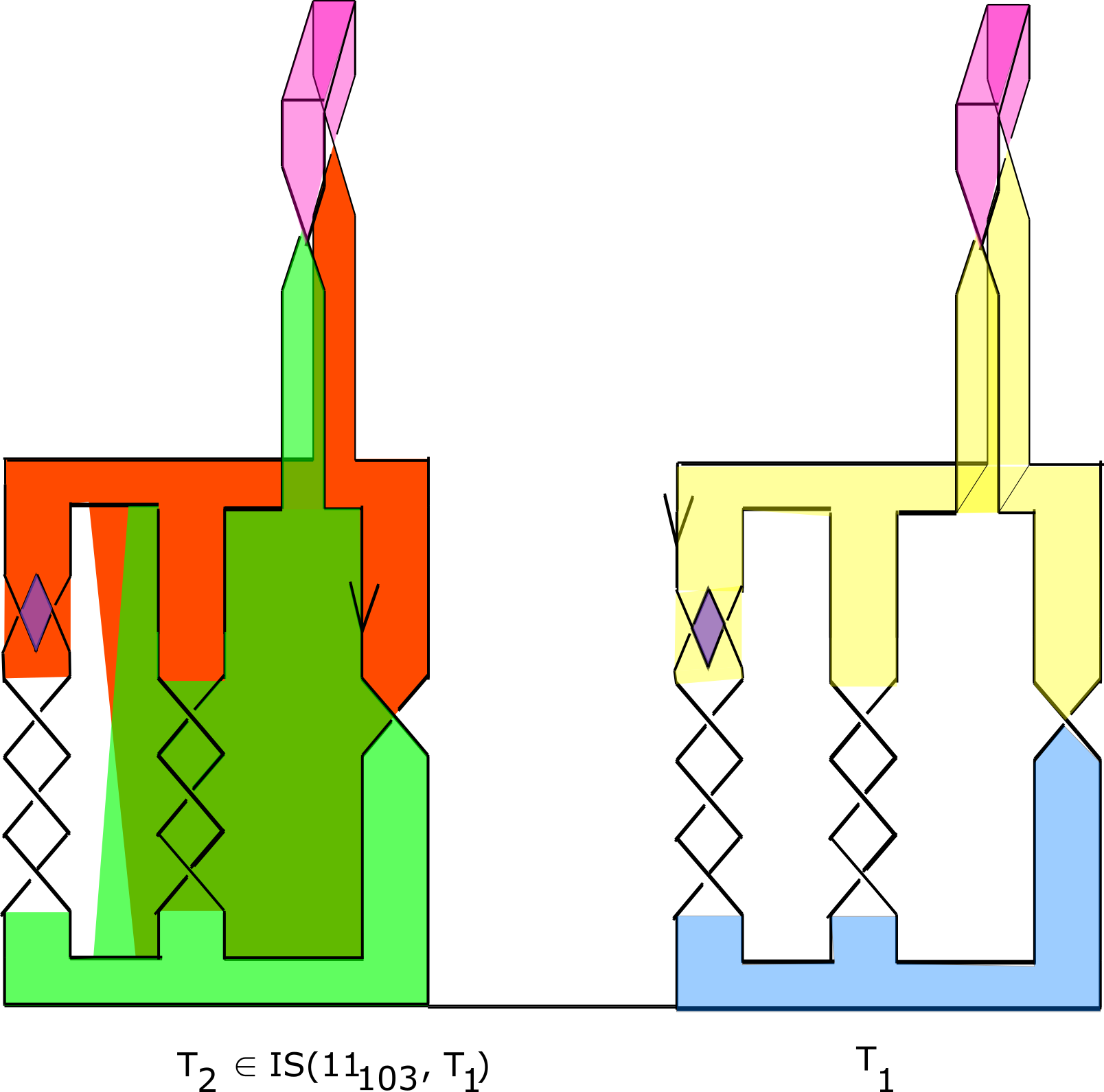}
  \caption{}
  \label{fig:boat1}
\end{figure}

\newpage

\section*{Finding $IS(L,T_2)$}

\begin{theorem}

\vspace{1cm}

Let $L = 11_{103},$ and let $D$ be an oriented reduced alternating  diagram of $L.$  Let $T_1$ be the surface obtained by applying Seifert's algorithm to $D.$ Let $IS(L,T_1) = T_2.$

Then  $$IS(L,T_2) = T_1.$$

\end{theorem}

From the previous Lemma, 
$$IS(L,T_1) = T_2$$ was obtained by applying 2 product decomposition to the surface $T_1.$

Let $H_1$ and $H_2$ be the Hopf bands plumbed to the surface $S_1(7_4).$ The surfaces $T_1$ and $T_2$ are parallel along $H_1.$ [picture]

$$L = 7_4\cup_{D_1}\partial(H_1) \cup_{D_2}\partial(H_2),$$ Consider the boundary of the surface $$\Tilde{L} = \partial(S_1(7_4)\cup_{D_2}(H_2))$$ 

From the previous lemma, we have:

$IS(L,T_1) = [T_2].$

If we deplumb $H_1$ from $T_i$, we denote the new surfaces (isotopy class) on $\tilde{L}$ as $S_i(\tilde{L}).$

$$\phi_1:IS(L,T_1)\to IS(\Tilde{L}, (S_1(7_4))\cup_{D_2} H_2) = IS(\Tilde{L}, (S_1(\widetilde{L}))$$

and the map $\phi_1$ is a bijection.  

$$[S_2(\Tilde{L})] : = \phi_1([T_2]).$$

and let $S_2(\Tilde{L})$ be a surface in the isotopy class  which is disjoint from $S_1(\Tilde{L}).$ 

We will use the  Lemma for

$$\Tilde{\phi}: IS(L,T_2)\to IS(\Tilde{L}, S_2(\tilde{L}))$$ is a bijection as well since we obtain $S_2(\widetilde{L})$  by deplumbing the Hopf band $H_1$ (equivalent to product decomposition on $H_1$)

\emph{We need to prove:}

$$IS(\Tilde{L},S_2(\Tilde{L})) = S_1(\Tilde{L}) = S_1(7_4)\cup_{D_1} H_2.$$

\begin{proof}

$ $

\vspace{1cm}

Let $[F] \in IS(\Tilde{L}, S_2(\Tilde{L}))$ and $[F]\neq [S_1(\Tilde{L})].$

Note that if $F\cap S_1(\Tilde{L}) = \phi$ this implies either 

\begin{itemize}
    \item $[F] = [S_1(\Tilde{L})]$ (which contradicts the assumption on $[F]\neq [S_1(\Tilde{L})]$); or
    
    \item  $[F]\in IS(\Tilde{L},S_1(\Tilde{L})) = [S_2(\Tilde{L})].$  (which contradicts the assumption that $[F]\in IS(\Tilde{L},S_2(\Tilde{L}))$)
    
    Hence $F$ intersects $S_1(\Tilde{L})$ essentially
    
    \end{itemize}

\begin{figure}[H]
  \includegraphics[width=7.5cm]{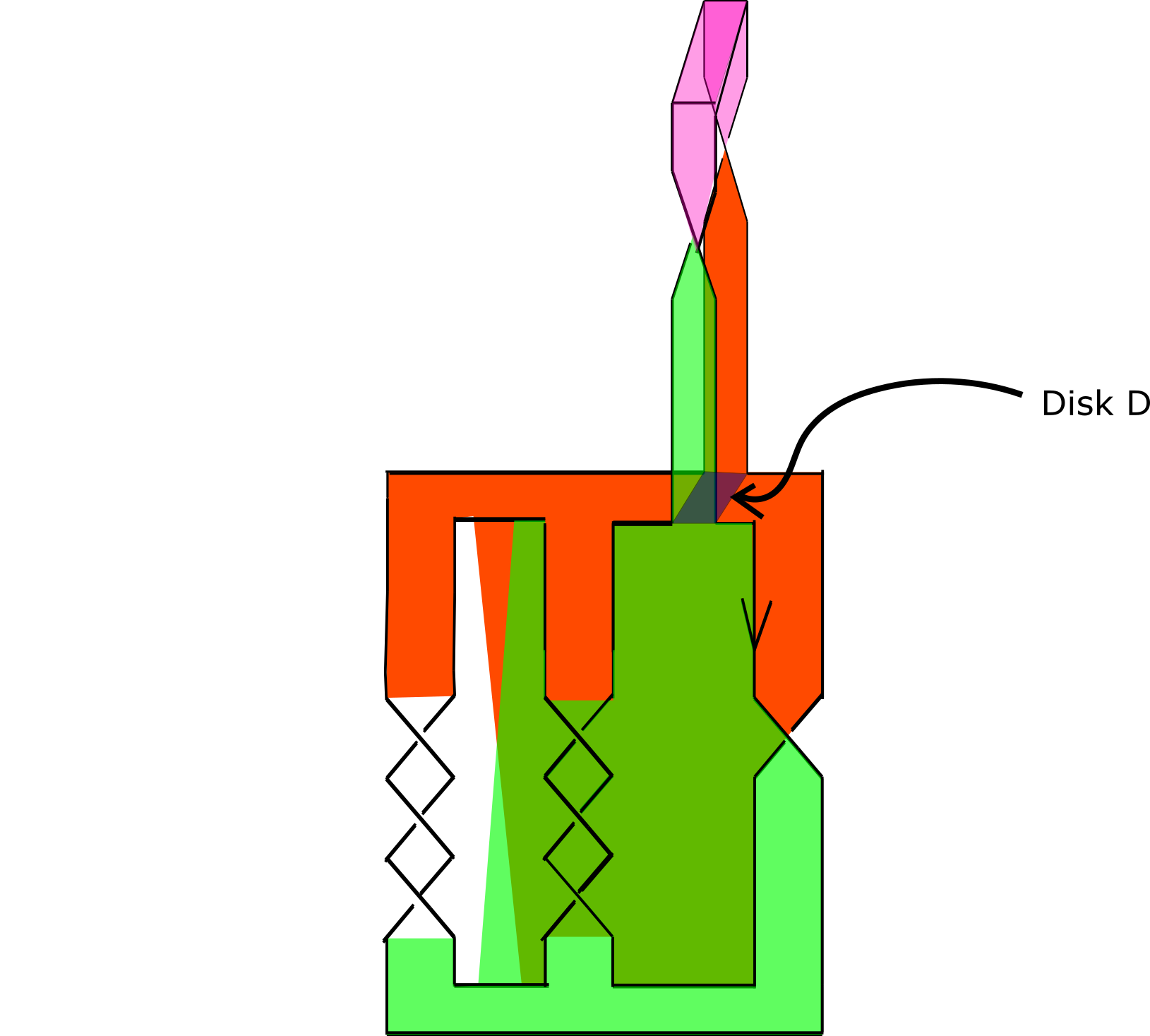}
  \caption{}
  \label{fig:boat1}
\end{figure}

Consider the two balls $V$ and $W$ such that $$V\cup W = S^3,$$ with $$V\cap (S_2(\Tilde{L})\cup D) = H_1 \text{ and  } W\cap (S_2(\Tilde{L})\cup D) = S_2(7_4),$$ where $D$ is 4-gon disk, boundary of which is a  rectangle with two of the opposite sides are $\partial V \cap S_2(\Tilde{L}),$ and  $S_2(\widetilde{L})$ is an embedded surface in $E(\widetilde{L}).$ 

\vspace{1cm}

[In the following pictures we assume the ball at the top is $V$ and the one at the bottom is $W.$]

\vspace{1cm}

Let $F$ be as above. $F \cap \partial V $ consists of two arcs since $F$ is a minimal genus Seifert surface and we assume that $F$ intersects $\partial V$ transversely. There are three cases as shown:  

\begin{figure}[H]
  \includegraphics[width=9cm]{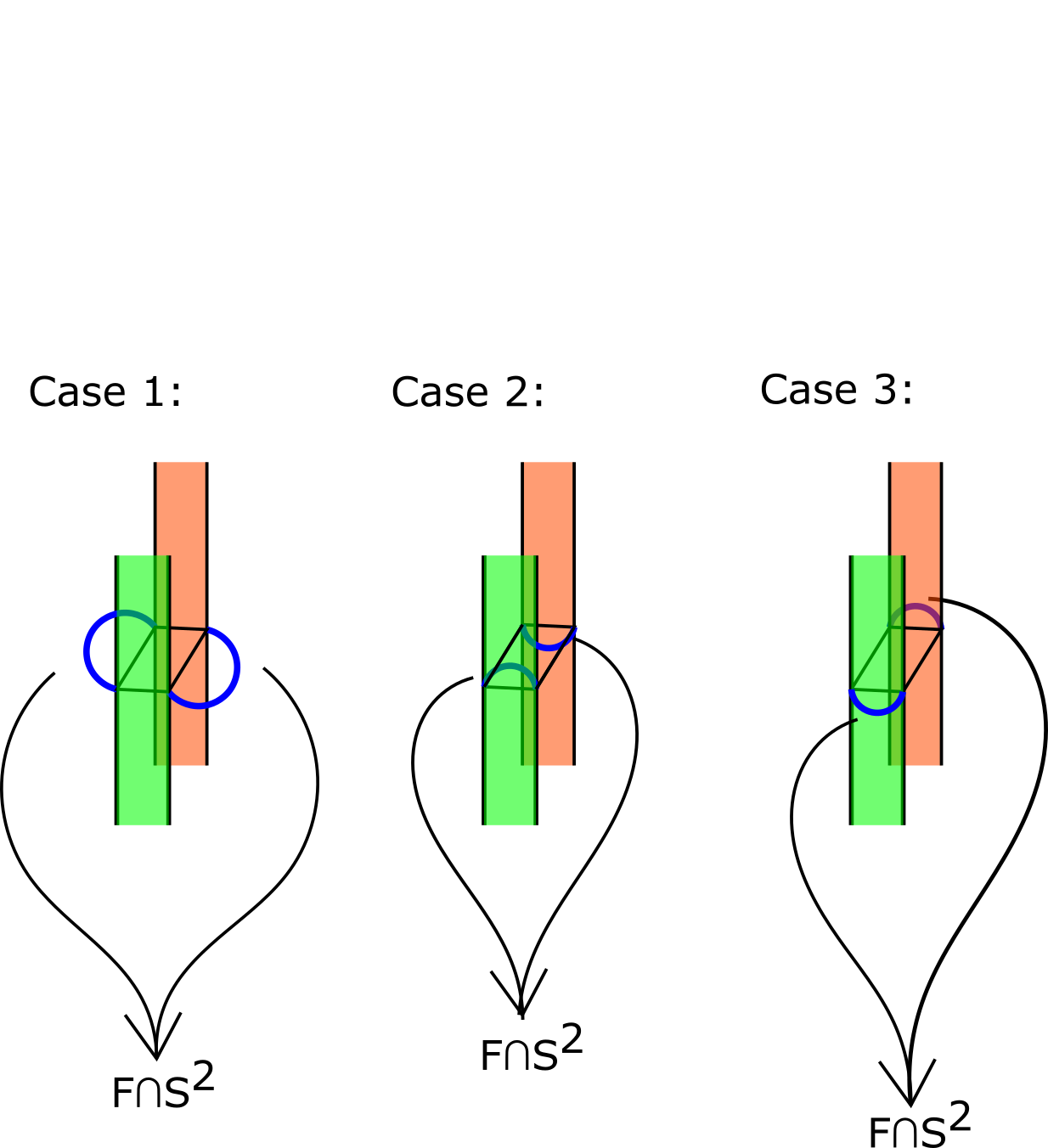}
  \caption{}
  \label{fig:boat1}
\end{figure}

\clearpage

We want to show that if $\Tilde{L}$ is the boundary of the surface $S_1(7_4) \cup_{D} H_2,$ then $IS(\Tilde{L},S_2(\tilde{L})) = S_1(7_4) \cup H_2$

\subsubsection*{Case 1:}

Consider $(M,\gamma)$ to be the complementary sutured manifold of $S_2(\widetilde{L}).$ Let $(M_1,\gamma_1)$ (resp. $(M_2,\gamma_2)$) be the complimentary sutured manifold for $(H_2,\partial H_2)$ (resp $(S_2(7_4) , 7_4)$). ($F,\partial F)\subset (M,\gamma)$ is a sutured surface in the complementary sutured manifold of $S_2(\widetilde{L}).$ Let $D_1$ and $D_1^c$ be the 4-gon's whose boundary is $((F\cap S^2)\cup (S_2(\tilde{L})\cap S^2))\cap (M,\gamma).$ [Let $D_1$ be the disk in the interior side of the diagram].

$$(\widetilde{M_1},\widetilde{\gamma_1}) = (M\cap V_1-N(D_1), (\gamma\cap V_1) \cup N(F\cap S^2)) \cong (M_1,\gamma_1).$$$$F_1 = F\cap V_1  \text{ is a } \widetilde{\gamma_1} \text{ surface in }\widetilde{M_1} \text{ along with } \partial F_1 \cong \partial H_2.$$

Since $H_2$ is a Hopf band, we know that $\partial H_2$ spans a unique Seifert surface and $ H_2$ is the fibred surface.

We isotope $F_1$ so that in $V_1$ it's parallel to $H_1.$ Let $e_1:[0,1]\times F_1\to \widetilde{M_1}$ be the isotopy with $e_1(0) = Id_{F_1} $ and $e_1(1) (F_1) = R_{(+)}(\widetilde{\gamma_1}),$ a parallel copy of $H_2.$ Since $H_2$ is fibred , we have an isotopy $e_2:[0,1]\times F_1\to \widetilde{M_1}$ with $e_1(0) = Id_{F_1} $ and $e_1(1) (F_1) = R_{(-)}(\widetilde{\gamma_1}),$ a parallel copy of $H_2.$

Note that $(\widetilde{M_1},\widetilde{\gamma_1})\cong (M_1,\gamma_1).$ $H_2$ is a unique fibred surface on $\partial H_2,$ and $F_1$ is a sutured surface on $(\widetilde{M_1}).$ This implies that $F_1$ is isotopic to $R_{+}(\widetilde{\gamma_1})$ or $R_{-}(\widetilde{\gamma_1}).$

Without loss of generality, let us assume that $F_1$ is isotopic to $R_{+}(\widetilde{\gamma_1}).$ Isotope the surface $F$ to $\widetilde{F}$ (consider an ambient sutured isotopy $e:I\times F\to M$ such that $e(0) = \text{Identity }$ and $e(1)(F) = \widetilde{F}$ such that $e(1)(F_1) = R_{+}(\gamma_1)$).

Hence we have $\widetilde{F}_1 = R_{+}(\widetilde{\gamma_1})$

Consider $I\times D_1 \subset \widetilde{M_1}$ with $\{0\}\times D_1 = D_1 \subset \partial V $ and the vertices of $\partial(\{t\} \times D_1) \in \widetilde{\gamma_1}$

$D_1' = \{1\}\times D_1 \in \widetilde{M_1}$ is an isotopic copy of $D$ in $\widetilde{M_1}.$

\vspace{1cm}

(If $F_1$ is isotopic to $R_{-}(\widetilde{\gamma_1}),$ then isotope the surface $F$ to $\widetilde{F}$ via an ambient sutured isotopy $e:I\times F\to M$ such that $e_0:F\to F \text{ is } \text{ the identity map }$ and $e_1(F) = \widetilde{F}$ with $e_1 (F_1) = R_{-}(\gamma_1)$).

We have $\widetilde{F}_1 = R_{-}(\widetilde{\gamma_1})$

Consider $I\times (D_1^{c} = \partial V - D_1^{\circ}) \subset \widetilde{M_1}$ with $\{0\}\times D_1^c = D_1^c \subset \partial V $ and the vertices of $\partial(\{t\} \times D_1^c) \in \widetilde{\gamma_1}$)

\vspace{1cm}

Let $V_1 = V - {(\partial V\times [0,1))}$ and $W_1 = W\cup {(\partial V\times [0,1])}$ with 

$$V_1\cap W_1 = \partial V_1 = \partial V \times \{1\}. [ D_1 \times \{t\}\subset \partial V \times \{t\}.]$$

Since $\widetilde{F}\cap V = \widetilde{F_1}$ is parallel to $H_2,$ $\widetilde{F}$ intersects the sphere $\partial V_1$ in two arcs perpendicular to $\widetilde{F_1}\cap V.$   

We have $[F] = [\widetilde{F}] \in IS(\Tilde{L}, S_2(\Tilde{L}))$which means $\widetilde{F}$ be a surface disjoint from $S_2(\Tilde{L}).$ and we know from the previous assumptions that $[\widetilde{F}]\neq [S_1(\Tilde{L})].$

Consider the two balls $V_1$ and $W_1.$ $$V_1\cap S_2(\Tilde{L})\cup D' = H_1 \text{ and  } W_1\cap S_2(\Tilde{L})\cup D' = S_2(7_4),$$ where $D = \{1\}\times D$ is the 4-gon disk, boundary of which is a  rectangle with two of the opposite sides are $\partial V_1 \cap S_2(\Tilde{L}).$

\vspace{.5cm} 

Since $\widetilde{F}$ intersects $D'$( or $D'^{c}$) in perpendicular arcs to $\widetilde{F}\cap D$( or $\widetilde{F}\cap D^{c}$) , Case 1 reduces to Case 2( or Case 3). 

\subsubsection*{Case 2 (Similarly for Case 3):}

$F\subset (M,\gamma)$ is a sutured surface in the complementary sutured manifold of $T_2.$ Case 2 implies that $F\cap S^2$ is a pair of arcs parallel to a copy of $T^2\cap S^2$ in $R_{+}(\gamma).$ [Case 3 indicates the same in $R_{-}(\gamma).$] 

$S^2$ separates $S^3$ in two balls $V_1$ and $V_2$ respectively [$V_1$ denotes the ball at the top in the diagram.]

Filling $V_1$ to a point, we isotope $S^3$ to $V_2,$ adjoining a point, a new $S^3.$ We have the link $L = 7_4$ and we have a Seifert surface $F_2 = F\cap V_2$  with boundary link $L = 7_4.$ Let $r_1$
and $r_2$ be 2 arcs in $S^2$ perpendicular to the arcs $F_2\cap S^2.$ Let the 4-gon  in $S^2$ bounded by $r_1,r_2, F\cap S^2$ be $D_2.$

(For Case 3: the 4-gon is bounded by the same arcs. $D_2^c$ satisfies the same condition, in the diagram it is the complement (in $S^2$) of the 4-gon in Case 2).

$$(\widetilde{M_2},\widetilde{\gamma_2}) = (M\cap V_2-N(D_2), (\gamma\cap V_2) \cup N(F\cap S^2)).$$

Let $(M_1,\gamma_1)$ (resp. $(M_2,\gamma_2)$) be the complementary sutured manifold for $(H_2,\partial H_2)$ (resp $(S_2(7_4) , 7_4)$). As mentioned before $F_2 = F\cap V_2$ is a $\widetilde{\gamma_2}$ surface in $\widetilde{M_2}$ along with $\partial F_2 = 7_4.$

Note that: $$(\widetilde{M_2},\widetilde{\gamma_2})\cong (M_2,\gamma_2).$$

This implies $F_2$ in $\widetilde{M_2}$ is isotopic to $S_1(7_4)$ or $S_2(7_4).$

\subsubsection*{Case a)}

$F_2$ in $\widetilde{M_2}$ is isotopic to $S_1(7_4).$ Let $S_1(\widetilde{L})$ be the standard embedded surface in $(M,\gamma)$ such that $S_1(\widetilde{L})$ restricted to $V_1$ is in $R_{+}({\gamma}).$  Let $S_1 = (S_1(\widetilde{L})\cap V_2)\cup D = S_1(7_4)$ (where $D$ is the plumbing disk) be the standard embedded surface in  $\widetilde{M_2}.$ Since $F_2$ is isotopic to $S_1(7_4)$ in $\widetilde{M_2}$ and both the surfaces are disjoint from $T_2$, we have an isotopy $e:I\times F_2\to \widetilde{M_2}$ such that $e(0,F_2) = F_2$ and $e(1,F_2)$ is a parallel copy of $ S_1(7_4)$ in $\widetilde{M_2}.$ 

\vspace{.3cm}

A priori we have the assumption that $F$ intersects $S_1(\widetilde{L})$ essentially. We would like to show that we can isotope $F$ in $E(L)$ such that $F$ is disjoint from $S_1(\widetilde{L}).$ That would prove the fact that $F$ is isotopic to $S_1(\widetilde{L}) = S_1(7_4)\cup_D H_2.$ where $\cup_D$ denotes the plumbing of the two surfaces. 

Let $S_1(\widetilde{L})$ be the embedded surface. As before, $V_1$ and $V_2$ denote the balls partitioning $S^3$ with $V_1$ being the ball at the top and $V_2$ at the bottom of the diagram with the plumbing disk $D\subset \partial V_1$. $F_2$ in $E(L)$ is isotopic to $S_1(7_4).$ Let $(M,\gamma)$ be the complementary sutured manifold of $S_1(\widetilde{L}).$ There is an isotopy $e:I\times F_2\to M_2$ such that $e_1(F_2)$ is parallel to $R_{+}(\gamma)$ or $R_-{(\gamma)}$ and $e_1(F_2)$ does not intersect $S_1(\widetilde{L}).$ Without loss of generality, we can assume $e_1(F_2) = R_+(\gamma).$ 

Let $N(S^2) = S^2\times[-1,1]$ with $S^2\times \{0\} = \partial V_1$ and $S^2\times \{-1\} \in V_2$ be a regular neighborhood of the Murasugi sphere $S^2.$

Let $V_{\text{down}} = V_2 - (N(S^2)^{\circ})$ and $V_{up} = \overline{V_1\cup N(S^2)}$ 
Let's consider 2 arcs $\alpha$ perpendicular to $S^2\times \{-1\} \cap e_1(F_2)$ in $S^2\times \{-1\}$. Let $R$ be a rectangular disk with the boundary being $\alpha$ and $S^2\times \{-1\} \cap e_1(F_2).$ $F_2' = (F\cap V_{down}) \cup R$ is a sutured surface in $(M_2' = V_{\text{down}}\cap M,\gamma_2')$ with the boundary link of $F_2' $ being $7_4.$ Then $F_2'$ is parallel to $R_+(\gamma_2')$  since $F_2$ is isotopic to $R_+(\gamma),$  isotope $F$ by an ambient isotopy such that $F_2' = R_+(\gamma_2').$

Consider $(M_1' = M\cap V_{\text{up}},\gamma_1')$ with $F_1' = F\cap V_{up}.$ being a sutured surface isotopic to the Hopf band . Hence we have an isotopy $e^{1'}:I\times F_1' \to E(L)$ such that $e'(1) $ is parallel to $R_+(\gamma_1).$ (since Hopf band is fibred). Isotope $F$ such that $F\cap(V_1\cup S^2\times [-1/2,1]) \subset R_+(\gamma)$ and by isotopy extension we can demand  $F_2' = R_+(\gamma_2').$ Hence, the surface $F$ intersects  $S_1(\widetilde{L})$ $\in V_2$ in $S^2\times (-1/2,-1).$ Since $F$ is a minimal genus Seifert surface, there is a product region between $F$ and $S_1(\widetilde{L}). $ Hence $F$ can be made disjoint from $S_1(\widetilde{L})$ which implies that $F$ is isotopic to $S_1(\widetilde{L}).$

\subsubsection*{Case b)}

$F_2$ is isotopic to $S_2(7_4).$ In this case we know that $F_2$ is parallel to $R_{+}(\widetilde{\gamma_2})$ or $R_{-}(\widetilde{\gamma_2}).$ WLOG let's assume $F_2$ is isotopic to $R_+(\widetilde{\gamma_2}).$ Hence there is an sutured isotopy $e:F_2\times I\to \widetilde{M_2}$ such that $e_0(F_2) = F_2$ and $e_1(F_2) = R_+(\widetilde{\gamma_2})).$ This implies that there is a product region between $F_2$ and $R_{+}(\widetilde{\gamma_2}).$ 

Consider a thin product, $S^2\times I\subset V_2$ such that $S_2\times \{0\} = \partial V_2$ and $S_2\times \{1\}\subset V_2.$

Let $$V_{up} = {V_1 \cup (S^2\times I)};\hspace{1cm} V_{down}  = \overline{V_2 - (S^2\times [0,1])}. $$

This implies 

$V_{up}$ and $ V_{down}$ are balls with 
$\partial V_{up} = \partial V_{down} = S^2\times \{1\}.$

We move $F$ by a sutured isotopy such that $\widetilde{M_1}$ is fixed and its restriction to $F_2$ is $e|_{F_2 \times [0,r]}$ , $r$ close to $1.$ 

Let $(M^*,\gamma^*) = (M\cup V_{\text{down}} , (\partial V_{up} - M) \cup \gamma \cap V_{up})).$ Consider $F^*$ be the sutured surface obtained by adding two rectangles along the arcs $F\cap \partial V_{up}.$ We know that $F^{*}$ is a sutured surface on $(M^*,\gamma^*)$ homeomorphic to $(M_1,\gamma_1.)$ Since Hopf band is the unique surface spanned by the core of $\gamma_1 = \text{Hopf Link}$, $F^*$ is isotopic to Hopf band, (let $e^*$ be the isotopy $e^*: F^*\times I \to M^*$). There are 2 possibilities $e^*_1(F^*) = R_+(\gamma^*)$ or $R_-(\gamma^*).$  Since the Hopf band is fibred,  we have an isotopy from $R_+(\gamma^*)$ to $R_-(\gamma^*).$ We consider the isotopy $e^*_1(F^*) = R_+(\gamma^*).$ [If apriori, $F^*$ is isotopic to $R_-(\gamma^*)$, then we can compose with the isotopy arising from the fibredness of the Hopf link]. 

So we have two embeddings, 

$$e((F\cap V_{\text{down}})\times [0,1])\subset M\cap V_{\text{down}} ; e^*((F\cap V_{up})\times [0,1]) \subset M\cap V_{up}$$ and 

$$e|_{F\cap \partial V_{up}\times [0,1]} = e^{*}|_{F\cap \partial V_{up}\times [0,1]}.$$

hence we can connect these two embeddings to get a resulting embedding $\widetilde{e}$ such that $\widetilde{e}_0 = id$ and $\widetilde{e}_(1) (F) = R_{(\gamma)},$ which proves that $F$ is isotopic to $\widetilde{S_2}$ contradicting the fact that $F\in IS(\widetilde{L}, \widetilde{S_2}).$

This concludes that $IS(\widetilde{L}, \widetilde{S_2}) = \widetilde{S_1}.$

and 
$$IS(L,T_2) = T_1.$$

Hence $IS(11_{103}) = \{T_1,T_2\}$

 \begin{tikzpicture}
    \draw[olive,thick,latex-latex] (0,0) -- (7,0)
    node[pos=0,mynode,fill=red,label=above:\textcolor{red}{$T_1$}]{}
    node[pos=1,mynode,fill=blue,text=blue,label=above:\textcolor{blue}{$T_2$}]{};
  \end{tikzpicture}

\end{proof}

\section*{Special Case $11_{201}$}

\begin{figure}[H]
  \includegraphics[width=16cm]{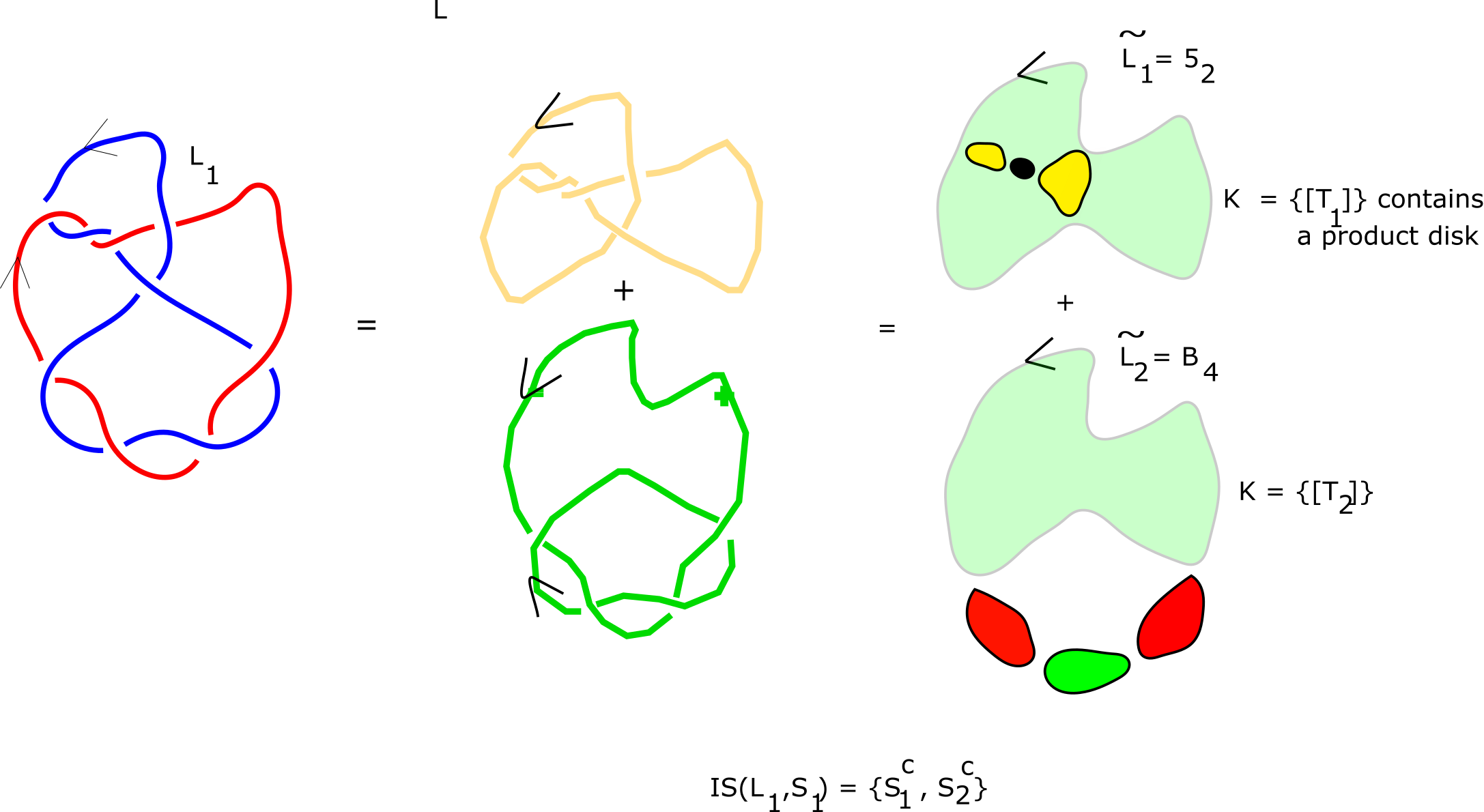}
  \caption{}
  \label{fig:boat1}
\end{figure}

\begin{figure}[H]
  \includegraphics[width=15cm]{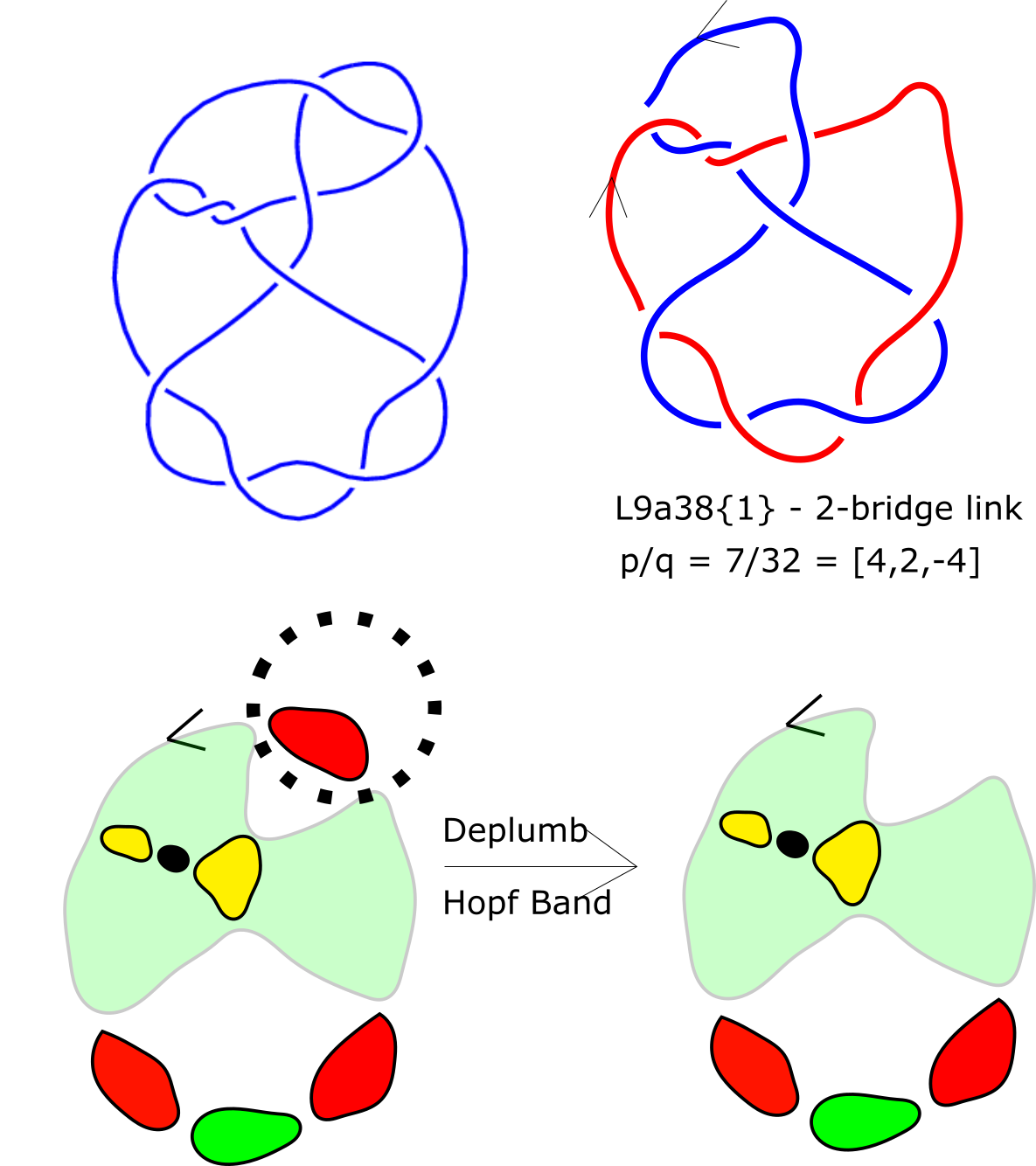}
  \caption{}
  \label{fig:boat1}
\end{figure}

Let $D$ be an oriented alternating diagram of $11_{201}$ and let $S$ be the surface obtained by applying Seifert's algorithm on $D.$ $S = S_1 \cup S_2$ where $S_2$ is the Hopf band. Since the Hopf band is a fibred surface on the Hopf link, hence we have, 

$IS(L,S)\cong IS(L_1,S_1).$

$S_1$ is a plumbing of $2$ surfaces

$S_1 = T_1\cup_{D} T_2,$ two unique Seifert surfaces on special alternating links $\tilde{L_1},$ and $\tilde{L_2}$ respectively. 

If $(M,\gamma)$ be the complementary sutured manifold for $S_1$ and $(M_1,\gamma_1,A_1)$ along with $(M_2,\gamma_2,A_2)$ be the marked complementary sutured manifolds for $T_1$ and $T_2$ respectively, then $M_2$ contains a product disk with $A_2$ as an edge. By the theorem of Kakimizu (the plumbed surface $S_1$ satisfies the rest of the condition and $S_1 = \tilde{S}_1$) we know, 

$$IS(L_1) =  \begin{tikzpicture}
    \draw[olive,thick,latex-latex] (0,0) -- (1,0)
    node[pos=0,mynode,fill=red,label=above:\textcolor{red}{$[S_1^c]$}]{}
    node[pos=1,mynode,fill=blue,text=blue,label=above:\textcolor{blue}{$[S_1]$}]{};

    \draw[olive,thick,latex-latex] (1,0) -- (2,0)
    node[pos=0,mynode,fill=blue,label=above:\textcolor{blue}{$[S_1]$}]{}
    node[pos=1,mynode,fill=green,text=green,label=above:\textcolor{green}{$[\tilde{S}_1^c]$}]{};
  \end{tikzpicture} $$

Therefore $IS(L_1,S_1) = \{[S_1^c] , [S_2^c]\}.$

$IS(L,S)\cong IS(L_1,S_1).$

This implies that $IS(L,S) = \{{R_1} , {R_2}\}.$ Via the congruence, one of the surface, say $R_1 = S_1^c \cup_D S_2.$ The other surface is a parallel surface, parallel near the plumbing disk with the Hopf band. [It is similar to the surface $T_2$ found in the previous case.]

Since $R_1 = S_1^c \cup_D S_2$ and $S_2$ is a Hopf band,  $IS(L,R_1) \cong IS(L_1, S_1^c) = \{S_1\}. $ Hence $IS(L,R_1) = S_1\cup S_2.$

$S_2$ is a Hopf band and $R_2$ is a parallel surface with $S_2$ and $L_1$ 

The same proof as finding $IS(L,T_2)$ in the previous section applies and we have 

$IS(L,R_2) = S.$

Hence the Kakimizu complex of $K = 11_{201}$ is:

\begin{tikzpicture}
    \draw[olive,thick,latex-latex] (0,0) -- (1,0)
    node[pos=0,mynode,fill=red,label=above:\textcolor{red}{$[R_1]$}]{}
    node[pos=1,mynode,fill=blue,text=blue,label=above:\textcolor{blue}{$[S]$}]{};

    \draw[olive,thick,latex-latex] (1,0) -- (2,0)
    node[pos=0,mynode,fill=blue,label=above:\textcolor{blue}{$[S]$}]{}
    node[pos=1,mynode,fill=green,text=green,label=above:\textcolor{green}{$[R_2]$}]{};
  \end{tikzpicture} 

\begin{figure}[H]
  \includegraphics[width=15cm]{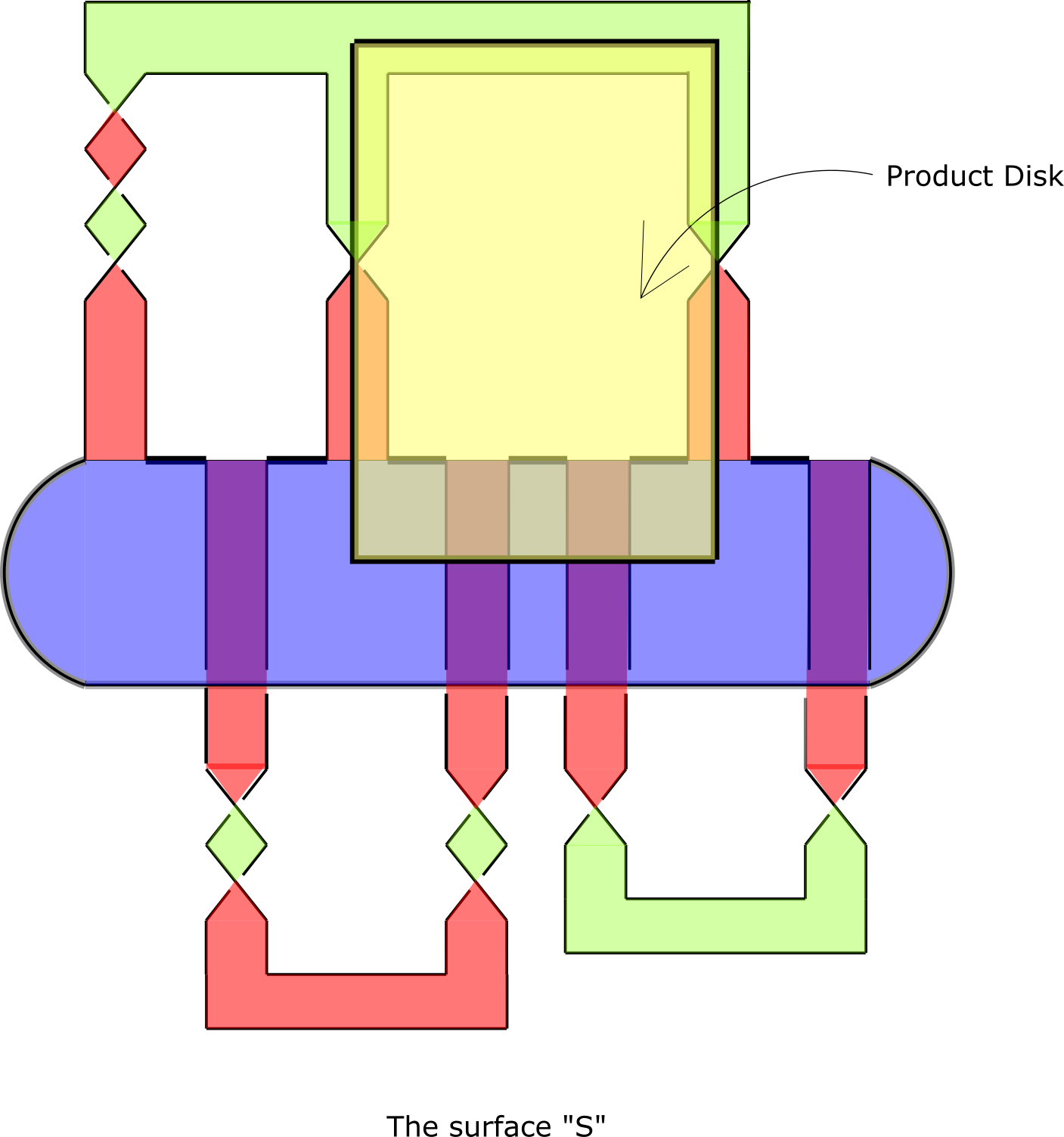}
  \caption{}
  \label{fig:boat1}
\end{figure}

\begin{figure}[H]
  \includegraphics[width=15cm]{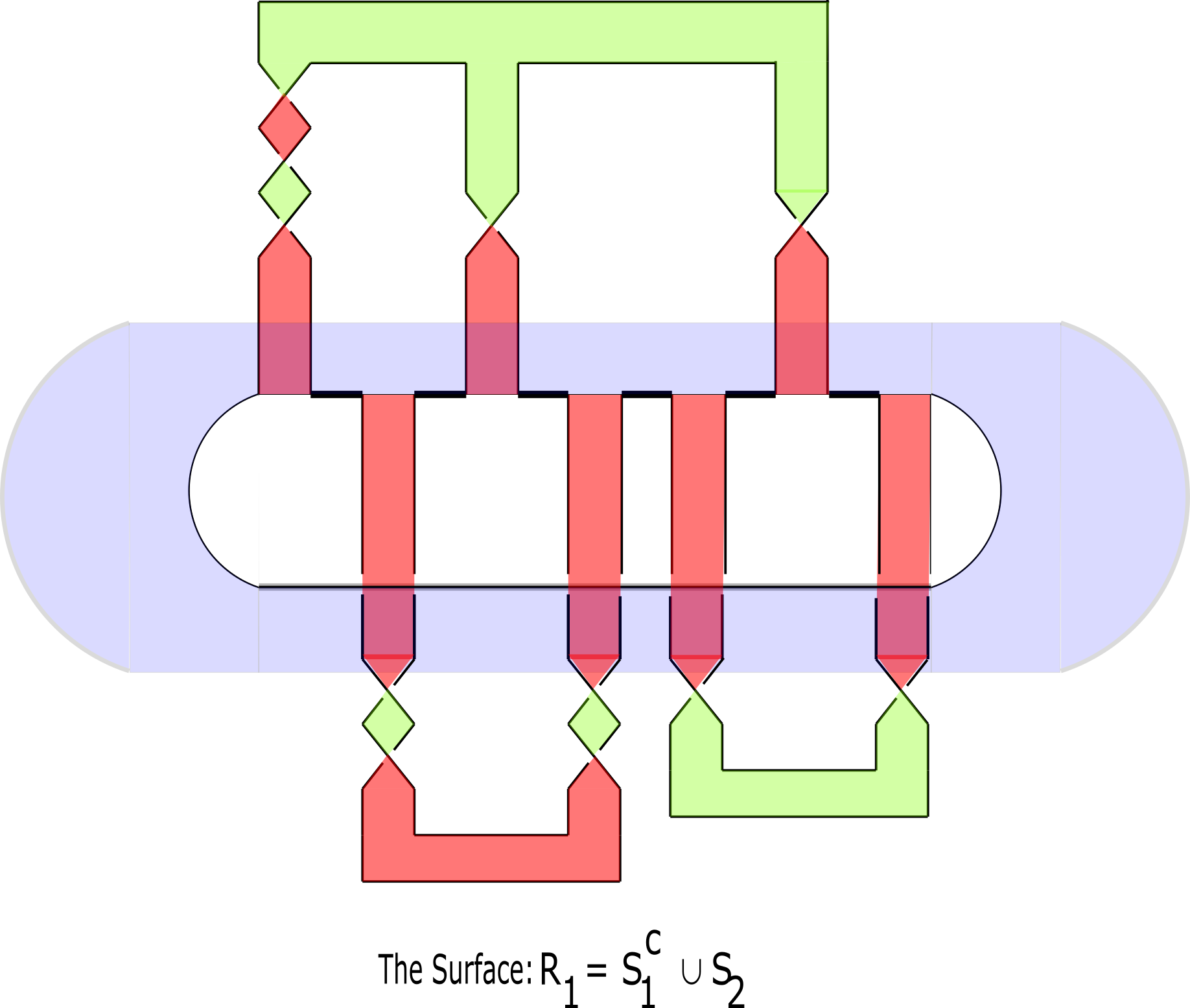}
  \caption{}
  \label{fig:boat1}
\end{figure}

\begin{figure}[H]
  \includegraphics[width=13cm]{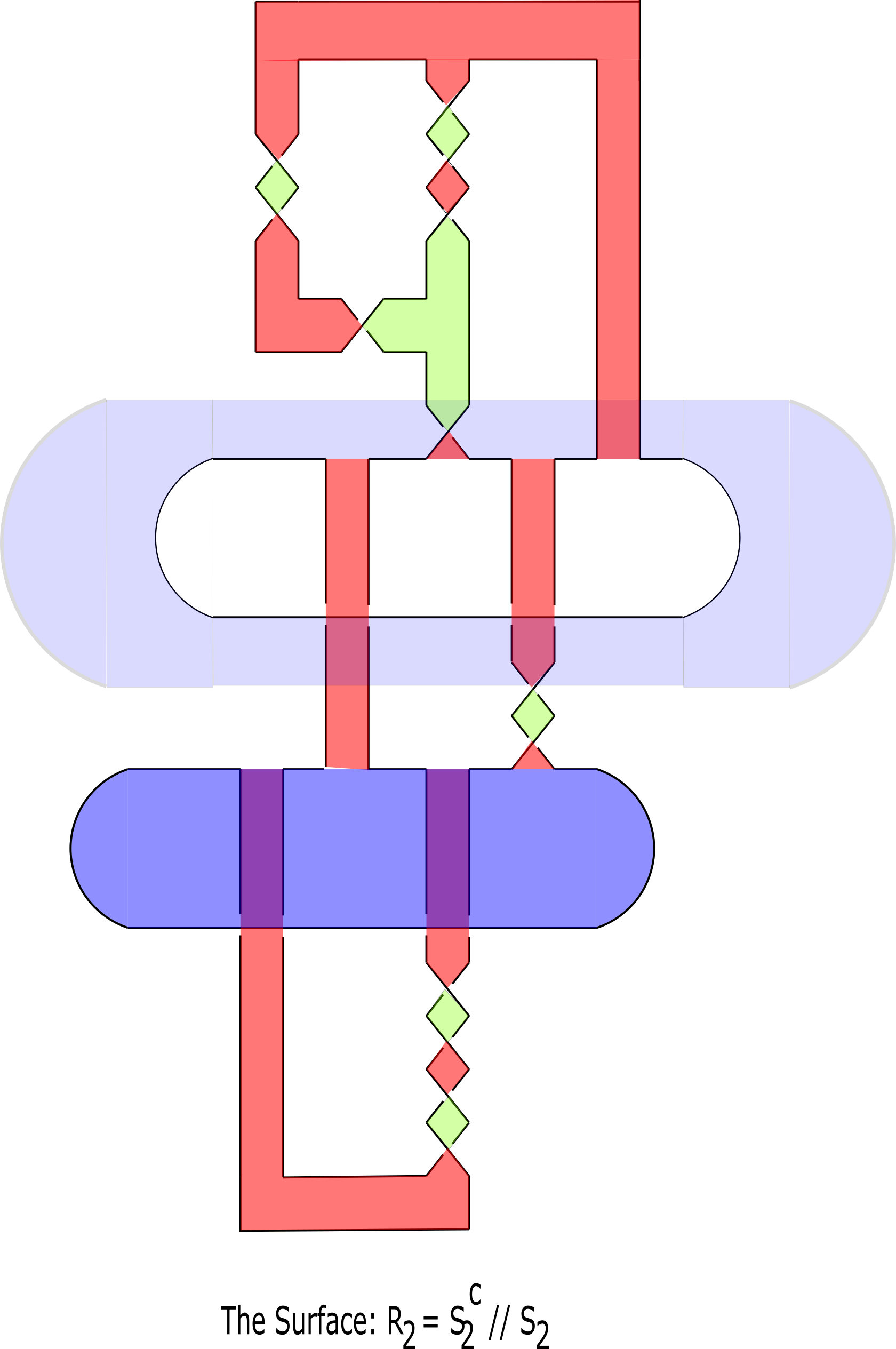}
  \caption{}
  \label{fig:boat1}
\end{figure}

\section{Further Questions}

The Kakimizu complexes of the class of $11$ crossing prime alternating knots is tabulated. 
\begin{itemize}
    \item The most immediate question is to find the data of Kakimizu complexes of 12 crossing prime, alternating knots to gain some structural information of the Kakimizu complex. Moreover, it seems plausible to extract from the current proofs, the Kakimizu complex of a special class of alternating links such that given an oriented diagram $D$ of $L,$ the Seifert surface obtained by applying Seifert's algorithm is a set of  plumbings of Seifert surfaces on special alternating links. Work is in progress on these 2 projects.

    \item The general goal is to find an algorithm to compute the Kakimizu complexes of every prime, non-split, oriented, alternating links.

    If possible, the goal is to extend to all hyperbolic links.

    \item Predict or provide more structure to the Kakimizu complex of a link. A general question is, given a graph $K,$ can $K$ be realized as a $1$ skeleton of the Kakimuzi complex of a link $L.$

\end{itemize}

\bibliography{refs}

\begin{thebibliography}{10}

\bibitem{10.32917/hmj/1150922486}
Osamu Kakimizu.
\newblock {Classification of the incompressible spanning surfaces for prime
  knots of 10 or less crossings}.
\newblock {\em Hiroshima Mathematical Journal}, 35(1):47 -- 92, 2005.

\bibitem{667e8b29-d581-32a3-ab92-1d9363b3b401}
Friedhelm Waldhausen.
\newblock On irreducible 3-manifolds which are sufficiently large.
\newblock {\em Annals of Mathematics}, 87(1):56--88, 1968.

\bibitem{WHITTEN1973373}
Wilbur Whitten.
\newblock Isotopy types of knot spanning surfaces.
\newblock {\em Topology}, 12(4):373--380, 1973.

\bibitem{knotinfo}
Charles Livingston and Allison~H. Moore.
\newblock Knotinfo: Table of knot invariants.
\newblock URL: url{knotinfo.math.indiana.edu}, Current Month Current Year.

\bibitem{https://doi.org/10.1112/jlms/s2-39.3.535}
P.~R. Cromwell.
\newblock Homogeneous links.
\newblock {\em Journal of the London Mathematical Society}, s2-39(3):535--552,
  1989.

\bibitem{Gabai1983TheMS}
David Gabai.
\newblock The murasugi sum is a natural geometric operation.
\newblock 1983.

\bibitem{banks2012minimal}
Jessica~E. Banks.
\newblock Minimal genus seifert surfaces for alternating links, 2012.

\bibitem{6129636b-b747-3d49-9531-0b1941964c84}
William Menasco and Morwen Thistlethwaite.
\newblock The classification of alternating links.
\newblock {\em Annals of Mathematics}, 138(1):113--171, 1993.

\bibitem{doi:10.1142/9789814529891}
Mikami Hirasawa and Makoto Sakuma.
\newblock {\em Minimal genus Seifert surfaces for alternating links. In KNOTS
  ’96 (Tokyo). World. Sci. Publ., River Edge, NJ, 1997.}, pages 383--394.

\bibitem{article}
Joshua Greene.
\newblock Alternating links and definite surfaces.
\newblock {\em Duke Mathematical Journal}, 166, 11 2015.

\bibitem{10.1112/blms/20.1.61}
Martin Scharlemann and Abigail Thompson.
\newblock {Finding Disjoint Seifert Surfaces}.
\newblock {\em Bulletin of the London Mathematical Society}, 20(1):61--64, 01
  1988.

\bibitem{25203298-73d8-3f83-8061-a5e5d88da654}
PIOTR PRZYTYCKI and JENNIFER SCHULTENS.
\newblock Contractibility of the kakimizu complex and symmetric seifert
  surfaces.
\newblock {\em Transactions of the American Mathematical Society},
  364(3):1489--1508, 2012.

\bibitem{Reidemeister1935HomotopieringeUL}
Kurt van Reidemeister.
\newblock Homotopieringe und linsenr{\"a}ume.
\newblock {\em Abhandlungen aus dem Mathematischen Seminar der Universit{\"a}t
  Hamburg}, 11:102--109, 1935.

\bibitem{Schubert1956}
Horst Schubert.
\newblock Knoten mit zwei brücken.
\newblock {\em Mathematische Zeitschrift}, 65:133--170, 1956.

\bibitem{Hatcher1985}
W.~Hatcher.~A, Thurston.
\newblock Incompressible surfaces in 2-bridge knot complements.
\newblock {\em Inventiones mathematicae}, 79:225--246, 1985.

\bibitem{Gabai1986}
David Gabai.
\newblock Detecting fibred links in s3.
\newblock {\em Commentarii mathematici Helvetici}, 61:519--555, 1986.

\end{thebibliography}
\bibliographystyle{unsrt}

\end{document}